\documentclass[a4paper,onecolumn,
superscriptaddress,10pt,nopdfoutputerror, accepted=2019-12-18,issue=1,
volume=2,
shorttitle=papers/compositionality-2-1]{compositionalityarticle}

\usepackage{hyperref}

\usepackage[T1]{fontenc}
\usepackage[english]{babel}

\usepackage{latexsym}
\usepackage{amsfonts}
\usepackage{amssymb}
\usepackage{amsmath}
\usepackage{amscd}
\usepackage{eucal}
\usepackage{amsthm}
\usepackage{pict2e}
\usepackage{pstricks}
\usepackage{pst-node}
\usepackage{verbatim} 
\usepackage{graphics}
\usepackage{graphicx}
\usepackage{color}   
\usepackage{framed}

\usepackage{bbm}

\psset{unit=1mm,arrowsize=2pt 3.5,arrowlength=1.1,linewidth=0.7pt,labelsep=2pt,
nodesep=2pt, arrowinset=0.7}
  
\def\1{\ensuremath{\mathbbm{1}}}%
\newcommand{\C}{\ensuremath{\mathbb C}} 

\renewcommand{\L}{\ensuremath{\mathbb L}}

\renewcommand{\P}{\ensuremath{\mathcal P}}

\newcommand{\D}{\ensuremath{\mathbb D}}
\newcommand{\E}{\ensuremath{\bb{E}}}

\newcommand{\F}{\ensuremath{\mathbb F}}

\newcommand{\Fop}{\ensuremath{\F\op}}

\newcommand{\cB}{{\cl{B}}}

\newcommand{\cE}{{\cl{E}}}
\newcommand{\cF}{{\cl{F}}}

\newcommand{\cV}{{\cl{V}}}

\newcommand{\bC}{\bb{C}}

\newcommand{\bI}{\bb{I}}

\newcommand{\bP}{\bb{P}}

\newcommand{\lt}{Lawvere theory}
\newcommand{\lts}{Lawvere theories}
\newcommand{\dl}{distributive law}
\newcommand{\dls}{distributive laws}

\newcommand{\Dls}{Distributive laws}
\newcommand{\fs}{factorisation system}
\newcommand{\fss}{factorisation systems}
\newcommand{\Fss}{Factorisation systems}

\newcommand{\iso}{\cong} 
\newcommand{\catequiv}{\simeq} 
\newcommand{\lra}{\ensuremath{\longrightarrow}} 

\newcommand{\dis}{\displaystyle}  

\renewcommand{\:}{\colon}

\renewcommand{\bf}{\bfseries}

\newcommand{\nd}{\hspace*{2pt}}

\newcommand{\fn}[1]{\ensuremath{\mbox{\bfseries {\upshape {#1}\hspace{1.3pt}}}}}

\newcommand{\cat}[1]{\ensuremath{\textrm{\bfseries {\upshape {#1}}}}}

\newcommand{\ob}{\mbox{\upshape{ob}\hspace{1.3pt}}}

\newcommand{\set}{\cat{Set}}
\newcommand{\Set}{{\cat{Set}}}
\newcommand{\Cat}{{\cat{Cat}}}
\newcommand{\Top}{{\cat{Top}}}

\newcommand{\Span}{{\cat{Span}}}
\newcommand{\Mod}{{\cat{Mod}}}
\newcommand{\Mnd}{{\cat{Mnd}}}
\newcommand{\Mon}{{\cat{Mon}}}
\newcommand{\Prof}{{\cat{Prof}}}
\newcommand{\CAT}{{\cat{CAT}}}

\newcommand{\PROF}{{\cat{PROF}}}
\newcommand{\ProfMon}{\Prof(\Mon)}
\newcommand{\ProfP}{\ensuremath{\Prof_{\cF}}}
\newcommand{\profp}{\ProfP}
\newcommand{\Profp}{\ProfP}
\newcommand{\finset}{\cat{FinSet}}
\newcommand{\FinSet}{\cat{FinSet}}

\newcommand{\Law}{\cat{Law}}
\newcommand{\Kl}{\fn{Kl}}

\newcommand{\scp}[1]{\scalebox{0.8}[0.7]{{\uppercase{#1}}}}

\newcommand{\dprof}{{\scp{prof}}}
\newcommand{\dfs}{{\scp{fs}}}
\newcommand{\dprofmon}{{\scp{profmon}}}
\newcommand{\dkl}{{\scp{kl}}}

\newcommand{\cl}[1]{\ensuremath{\mathcal {#1}}}
\newcommand{\bb}[1]{\ensuremath{\mathbb {#1}}}
\newcommand{\ed}{\end{document}}
\newcommand{\map}[1]{\ensuremath{\stackrel{{#1}}{\lra}}}

\newcommand{\demph}[1]{{\bfseries #1}}

\newcommand{\op}{^{\hspace{1pt}\textrm{\upshape{op}}}}

\newcommand{\uscore}{\ensuremath{\underline{\hspace{0.6em}}}}

\newcommand{\vv}[1]{\vspace*{#1}}
\newcommand{\hh}[1]{\hspace*{#1}}

\newcommand{\stra}{\psset{unit=0.1cm,nodesep=0pt,arrowsize=2pt 3.5,arrowlength=1.1} \pspicture(6.4,0) 
\pcline[linewidth=0.6pt]{->}(0.8,1.1)(5.8,1.1) \endpspicture}

\newcommand{\tra}{\psset{unit=0.1cm,nodesep=0pt,arrowsize=2pt 3.5,arrowlength=1.1} \pspicture(8,0)
\pcline[linewidth=0.6pt]{->}(1,1.1)(7,1.1) \endpspicture}

\newcommand{\tmapsto}{\psset{unit=0.1cm,nodesep=0pt,arrowsize=2pt 3.5,arrowlength=1.1} \pspicture(8,0) 
\pcline[linewidth=0.6pt]{|->}(1,1.2)(7,1.2) \endpspicture}

\newcommand{\tramap}[1]{\psset{unit=0.1cm,nodesep=0pt,labelsep=2pt,arrowsize=2pt 3.5,arrowlength=1.1} \pspicture(8,4)
\pcline[linewidth=0.6pt]{->}(1,1.1)(7,1.1)\naput{\ensuremath{\scriptstyle{#1}}} \endpspicture}

\newcommand{\tmap}{\tramap}

\newcommand{\ltramap}[1]{\psset{unit=0.1cm,nodesep=0pt,arrowsize=2pt 3.5,arrowlength=1.1} \pspicture(17,4)
\pcline[linewidth=0.6pt]{->}(2.5,1.1)(14.5,1.1)\naput{\ensuremath{\scriptstyle{#1}}} \endpspicture}

\newcommand{\ltmap}{\ltramap}

\newcommand{\prof}{\psset{unit=0.1cm,nodesep=0pt,arrowsize=2pt 3.5,arrowlength=1.1} \pspicture(10,0)
\pcline[linewidth=0.6pt,labelsep=2pt]{->}(1,1.1)(9,1.1)
\pcline[linewidth=0.6pt,tbarsize=3pt 1]{|-}(4.4,1.1)(8.5,1.1) \endpspicture}

\newcommand{\pra}{\prof}

\newcommand{\pmap}[1]{\psset{unit=0.1cm,nodesep=0pt,arrowsize=2pt 3.5,arrowlength=1.1} \pspicture(10,4)
\pcline[linewidth=0.6pt,labelsep=2pt]{->}(1,1.1)(9,1.1)\naput{\ensuremath{\scriptstyle{#1}}}
\pcline[linewidth=0.6pt,tbarsize=3pt 1]{|-}(4.4,1.1)(8.5,1.1) \endpspicture}

\newcommand{\dra}{\psset{unit=0.1cm,nodesep=0pt,arrowsize=2pt 3.5,arrowlength=1.1} \pspicture(10,0)
\pcline[linewidth=0.6pt]{->}(1,1.1)(9,1.1) \rput(5,1.1){\pscircle*(0,0){1.6pt}}
 \endpspicture}
 
\newcommand{\dmap}[1]{\psset{unit=0.1cm,nodesep=0pt,arrowsize=2pt 3.5,arrowlength=1.1} \pspicture(10,4)
\pcline[linewidth=0.6pt]{->}(1,1.1)(9,1.1)\naput{\ensuremath{\scriptstyle{#1}}} \rput(5,1.1){\pscircle*(0,0){1.6pt}}
 \endpspicture}

\newcommand{\cra}{\psset{unit=0.1cm,nodesep=0pt,arrowsize=2pt 3.5,arrowlength=1.1} \pspicture(10,0)
\pcline[linewidth=0.6pt]{->}(1,1.1)(9,1.1) \rput(5,1.1){\pscircle(0,0){1.6pt}}
 \endpspicture}
 
\newcommand{\cmap}[1]{\psset{unit=0.1cm,nodesep=0pt,arrowsize=2pt 3.5,arrowlength=1.1} \pspicture(10,4)
\pcline[linewidth=0.6pt]{->}(1,1.1)(9,1.1)\naput{\ensuremath{\scriptstyle{#1}}} \rput(5,1.1){\pscircle[linewidth=0.6pt](0,0){1.6pt}}
 \endpspicture}

\newcommand{\trta}{\psset{unit=0.1cm,nodesep=0pt} \pspicture(8,0)
\pcline[linewidth=0.6pt,doubleline=true,arrowinset=0.6,arrowlength=0.8,arrowsize=0.5pt 2.1]{->}(1,1.1)(7,1.1) \endpspicture}

\newcommand{\Tra}{\trta}

\newcommand{\noi}{\noindent}
\newcommand{\dend}{\ensuremath{\hfill \maltese}\end{mydefinition}}

\newcommand{\igc}[1]{\begin{tabular}[c]{c}\includegraphics[width={#1}]}
\newcommand{\igt}[1]{\begin{tabular}[t]{c}\\[-24pt] \includegraphics[width={#1}]}
\newcommand{\ei}{\end{tabular}}

\newtheorem{theorem}{Theorem}[section]

\newtheorem{proposition}[theorem]{Proposition}
\newtheorem{corollary}[theorem]{Corollary}

\newtheorem{definition}[theorem]{Definition}

\newtheorem{example}[theorem]{Example}
\newtheorem{note}[theorem]{Note}
\newtheorem{nonexample}[theorem]{Non-example}
\newtheorem{examples}[theorem]{Examples}
\newtheorem{remarks}[theorem]{Remarks}
\newtheorem{exercise}[theorem]{Exercise}
\newtheorem{question}[theorem]{Question}
\newtheorem{questions}[theorem]{Questions}
\newtheorem{algorithm}[theorem]{Algorithm}
\newtheorem{method}[theorem]{Method}

\newtheorem{defn}[theorem]{Definition}

\newtheorem{thm}[theorem]{Theorem}

\newtheorem{cor}[theorem]{Corollary}

\newenvironment{mydefinition}{\begin{definition}} {\end{definition}}
\newenvironment{myremark}{\begin{remark}\upshape} {\end{remark}}

\newenvironment{myremarks}{\begin{remarks}\upshape} {\end{remarks}}
\newenvironment{myexample}{\begin{example}\upshape} {\end{example}}

\newenvironment{myexamples}{\begin{examples}\upshape} {\end{examples}}

\newenvironment{myproposition}{\begin{proposition}\upshape} {\end{proposition}}

{ \end{sf}\end{framed}\end{minipage}
\end{center}}

\newtheoremstyle{example}{\topsep}{\topsep}%
     {}
     {}
     {\bfseries}
     {.}
     {8pt}
     {\thmname{#1}\thmnumber{ #2}\thmnote{ #3}}

   \theoremstyle{example}

\newtheorem{remark}[theorem]{Remark}
\newtheorem{eg}[theorem]{Example}

\newenvironment{prf}{\vspace{1ex}\begin{sloppypar}{\noindent\upshape
{\bfseries Proof.}}} {{\hspace*{\fill}
$\Box$}\end{sloppypar}\vspace{2ex}}

\newcommand{\glob}{
\xy
(-5,0)*+{.}="1";
(5,0)*+{.}="2";
{\ar@/^1pc/^{} "1";"2"};
{\ar@/_1pc/_{} "1";"2"};
{\ar@{=>}^{} (0,2)*{};(0,-2)*{}} ;
\endxy}

\newcommand{\numarabic}{\renewcommand{\labelenumi}{\arabic{enumi}.}}

\newcommand{\psinv}[6]{\xy
(#1,0)*+{#3}="x";
(#2,0)*+{#5}="y";
{\ar@<.7ex>^{#4} "x"; "y"};
{\ar@<.7ex>^{#6} "y"; "x"};
\endxy
}

\begin{document}

\title{Distributive laws for Lawvere theories}
\date{}
\author{Eugenia Cheng}
\email{echeng4@saic.edu}
\homepage{www.eugeniacheng.com}
\affiliation{School of the Art Institute of Chicago, Chicago, USA}

%
%


\maketitle

\begin{abstract}
Distributive laws give a way of combining two algebraic structures expressed as monads; in this paper we propose a theory of distributive laws for combining algebraic structures expressed as Lawvere theories.  We propose four approaches, involving profunctors, monoidal profunctors, an extension of the free finite-product category 2-monad from \Cat\ to \Prof, and factorisation systems respectively.  We exhibit comparison functors between {\bfseries CAT} and each of these new frameworks to show that the distributive laws between the Lawvere theories correspond in a suitable way to distributive laws between their associated finitary monads.  The different but equivalent formulations then provide, between them, a framework conducive to generalisation, but also an explicit description of the composite theories arising from distributive laws.
\end{abstract}

\tableofcontents


\section*{Introduction}
\addcontentsline{toc}{section}{Introduction}

\lts\ were introduced in \cite{law2} and were a great breakthrough in the
understanding of algebraic theories.  They give a different viewpoint from that 
of monads in how they implement the notion of arity.  One practical 
advantage of \lts\ over monads highlighted in \cite{hp1} is that 
\lts\ allow us to study models in 
different categories, starting from the same \lt.  For example, topological 
groups and ordinary groups both arise as models for the \lt\ for groups, whereas 
using monads we have to construct a monad on \Set\ for groups and a different 
(albeit related) monad on \Top\ for topological groups.

\Dls\ give us a way of combining algebraic theories expressed as monads.  The 
classic example combines the monad for Abelian groups and the monad for monoids 
(both monads being on \Set) to yield the monad for rings as the ``composite'' 
algebraic theory: the \dl\ makes the composite of the two monads into a new 
monad.  The theory for combining three or more monads is developed in 
\cite{che17}.  

It is well-known that \lts\ and monads are related---\lts\ correspond to \emph{finitary} monads on \Set.  This should not be thought of as a statement that \lts\ are ``merely'' a special case of monads; the above comments about models shows one way in which \lts\ are of importance in their own right.

A natural question then arises---is there a notion of \dl\ for \lts? Of course, given the above correspondence with finitary monads on \Set, one could simply say ``a \dl\ for \lts\ is a \dl\ between the associated finitary monads on \Set.''

However, we seek a formulation that is ``native'' to the framework of \lts.  In 
this paper we will provide four equivalent formulations at varying levels of 
abstraction.  As usual we expect the most abstract one to be more useful for 
theorising, and expect the most concrete one to be more useful for 
applications. Thus their equivalence should not be taken to mean that any of 
the definitions is redundant.

Our three most abstract formulations will come from observing that \lts\ may themselves be thought of as monads inside some other bicategory.  Having expressed \lts\ in this way it is natural to define \dls\ for \lts\ as \dls\ between the monads in these bicategories.  The bicategories in question are
\begin{enumerate}
 \item \Prof---categories, profunctors and natural transformations.
\item \ProfMon---as above but internal to monoids.
\item \ProfP---the Kleisli bicategory for the free finite-product category 2-monad extended from \Cat\ to \Prof.
\end{enumerate}
The advantage of (1) is that the bicategory \Prof\ is well-known and quite easy to understand; however not \emph{all} monads in here are \lts\, even if we restrict to the correct underlying 0-cell.

The approach using (2) is in some ways more naturally-arising than (1) and in 
fact helps us understand it.  Also, it is closely related to Lack's work on 
\dls\ for PROPs \cite{lac4}. 

The advantage of (3) is that, once we restrict to the correct underlying 0-cell, \emph{all} monads are \lts. It is this that enables us to prove that the composite monad in each of these three frameworks is also a \lt---it is immediate in (3) and then by the equivalence of the three definitions, the result will follow for (1) and (2).

Another advantage of (3) is that although (or because) this bicategory is very much harder to work with, it affords not only the most precision but also greater flexibility.   We will see that monads on other 0-cells may be thought of as ``typed'' \lts, and the setting also opens the possibility for changing the 2-monad $\cF$ to study different types of theory; this insight is all gained from Hyland \cite{hyl1}.

For the most concrete formulation, we unravel (1) and express it in terms of \fss.  The notions are equivalent, but the framework feels quite different from the above abstractions and therefore provides different insights.  For example, \dls\ for monads seem suited to considering composition of monads, whereas \fss\ seem suited to considering \emph{decompositions}.

Note that it is quite easy to make a wrong definition of \dl\ for \lts\ along the above lines, by working in an ill-chosen bicategory.  For example, every \lt\ is a monad in \Span\ (which is, after all, related to \Prof), but considering \dls\ in this bicategory gives the wrong notion, as we will show in Section~\ref{factsystems}.

As evidence that our definitions do give the correct notion, we prove that all our definitions of \dl\ for \lts\ correspond suitably to \dls\ between the associated finitary monads, with the composite \lts\ corresponding to the composite monads.  \emph{En passant}, we shed some more, abstract, light on the monad/\lt\ correspondence.  

Note that the tensor product of Lawvere theories is a way of combining Lawvere theories that is different from distributive laws.  The tensor product of two Lawvere theories always exists, whereas there is not always a distributive law of a given Lawvere theory over another.  It is said that in the tensor product ``all the operations of one theory commute with all the operations of the other'' but this must be understood in a particular sense: given an $m$-ary operation $f$ of the first theory and an $n$-ary operation $g$ of the second, in the tensor product $n$ copies of $f$ followed by $g$ is the same as $m$ copies of $g$ followed by $f$, as $(nm)$-ary operations.  This neither implies nor is implied by a distributive law.  For example, the theory of rings is not the tensor product of the theory of Abelian groups and the theory of monoids; the monad for rings is the composite of the monad for Abelian groups and the monad for monoids.  While this can be thought of as a type of commutativity between the group operation and the monoid operations, this is in a very different sense from the type of commutativity in the tensor product of Lawvere theories.

The paper is organised as follows.  In Section~\ref{lts} we briefly recall the definition of \lt\ and the correspondence with finitary monads on \Set.  In Section~\ref{dls} we briefly recall the notion of \dl\ between monads inside a bicategory.  Experts can skip both these sections with impunity.  In Sections~\ref{prof}--\ref{kleisli} we present our four different approaches to \dls\ for \lts\, and in Section~\ref{comparison} we provide the comparison.  We finish in Section~\ref{future} with some brief comments about the possibilities for future work.


\section{Lawvere theories}\label{lts}

In this section we recall the basic definitions and results about \lts\ that we will need in the rest of this paper.  Nothing in this section is new.  \lts\ were introduced in \cite{law2}; we find that \cite{hp1} gives a useful exposition.

The idea of a \lt\ is to encapsulate an algebraic theory as a category \L\ where
\begin{itemize}
 \item the objects of \L\ are the natural numbers, the ``arities'',
\item a morphism $k \tra 1$ is an operation of arity $k$, and
\item a morphism $k \tra m$ is $m$ operations of arity $k$.
\end{itemize}

\noi Let \F\ denote a skeleton of \finset, the category of finite sets and all functions between them.  So in particular the objects of \F\ are the natural numbers (including 0).

\begin{definition}
 A \emph{\lt} is a small category \L\ with (necessarily strictly associative) finite products, equipped with a strict product-preserving identity-on-objects functor
\[\alpha_\L \: \F\op \tra \L.\]
A \demph{morphism of \lts} from $\L$ to $\L'$ is a functor making the obvious triangle commute; note that such a functor necessarily strictly preserves finite products. Lawvere theories and their morphisms form a category \Law.

\end{definition}


\begin{myremark}
 It is worth making the structure of \F\ a little more explicit here as we will rely on this heavily later, especially when we consider the free finite-product category monad \cF \ in Section~\ref{kleisli}.  Since $\finset$ is equivalent to the free finite coproduct category on 1, $\F\op$ is equivalent to the free finite product category on 1.  Finite products are given by \emph{addition} of natural numbers, and so a morphism
\[\alpha \: k \tra m \in \F\]
is given by, for each $i \in [m]$, a choice of projection $k \tra 1$.  Hence $\alpha$ is precisely a function $[m] \tra [k]$ where we write $[k]$ for a set of $k$ elements. (We will sometimes omit the square brackets if confusion is unlikely.)

The idea for \lts\ is that \Fop\ encapsulates the operations that must generically exist in any algebraic theory: forgetting and repeating variables.  For each $m \in \F\op$ we have:
\begin{itemize}
 \item the $i$th product projection $m \tra 1$ corresponding to forgetting all $m$ variables except the $i$th one, and
\item the diagonal $1 \tra m$ corresponding to repeating a variable $m$ times.
\end{itemize}

\begin{mydefinition}
 The morphisms in \L\ are called \demph{operations}.
\end{mydefinition}

\end{myremark}

\begin{myexample}\label{ex1}
 
In the \lt\ for monoids, the 2-ary operations, that is, morphisms $2 \tra 1$, include the operations
\[ab, a, a^2, b, b^2, aba, ab^3a^5, \ldots \]
that is, everything in the free monoid on a 2-element set.  This could be seen as a different notion of arity from the one used to express algebraic theories via operads---in the (non-symmetric) operad for monoids the only 2-ary operation is $ab$; it could also be seen as a different notion of operation.  

A morphism $3 \tra 2$ is given by two 3-ary operations, eg
\[ \{ abc, ab^2c^2 \},  \{ bc^2a, ababc\}, \ldots \]
A typical composite looks like
\[3 \ltmap{\{abc, ab^2c^2\}} 2 \ltmap{\{x^2y\}} 1\]
yielding the composite 3-ary operation $abc\cdot abc\cdot ab^2c^2$.

Note that as a result of forgetting variables we have many different possible arities for the ``same'' operation.  For example starting with a 3-ary operation $abc$, say, we may precompose with variable-forgetting morphisms to express $abc$ as a $k$-ary operation where all variables apart from $a,b,c$ are forgotten: 
\[
\psset{unit=0.1cm,labelsep=2pt,nodesep=3pt}
\pspicture(0,0)(40,15)


\rput(55,0){\rnode{a1}{$1$}} 
\rput(55,13){\rnode{a3}{$3$}} 
\rput(40,13){\rnode{a4}{$4$}} 
\rput(25,13){\rnode{a5}{$5$}} 
\rput(10,13){\rnode{a6}{${\textcolor{white} 6}$}}
\rput(6,13){$\cdots$}

\ncline{->}{a3}{a1} \naput[npos=0.4]{{\scriptsize $abc$}}
\ncline{->}{a4}{a1}
\ncline{->}{a5}{a1}

\ncline[linestyle=dashed, dash=3pt 2pt]{->}{a4}{a3}
\ncline[linestyle=dashed, dash=3pt 2pt]{->}{a5}{a4}
\ncline[linestyle=dashed, dash=3pt 2pt]{->}{a6}{a5}

\endpspicture\]
As a 5-ary operation, for example, this might take the variables $a,b,c,d,e$ and return the operation $abc$.

\end{myexample}

\begin{myremark}
 There are many variations and generalisations of the notion \lt. Here are some examples. 

\begin{enumerate}

\item An alternate definition says that a \lt\ is any category with finite products in which every object is isomorphic to a finite cartesian power $x^n$ of a generic object $x$; this is then invariant under equivalence of categories.

\item Many-sorted theories: writing $\cF$ for the 2-monad 
for strictly associative products on the 2-category 
\Cat\ of small categories, and observing that $\cF 1 \catequiv \F\op$, we could 
instead use $\cF A$ for non-terminal categories $A$ to get \lts\ with sorts 
given by $A$.

\item Sometimes Lawvere Theories are known as finite-product theories, but sometimes finite-product theory is used to mean \emph{any} small category $C$ with finite products. In fact this notion of finite-product theory can be regarded as a special case of many-sorted theories in which the sorts are given by the objects of $C$.

\item Enriched theories: we could use enriched categories, and get a notion of enriched \lt, and higher-dimensional \lt; see \cite{pow3}.

\item $\Phi$-theories: we could use some other class $\Phi$ of limits than finite products, such as small products or finite limits; see \cite{lr1}.

\end{enumerate}

\end{myremark}

While \lts\ enable us to study, say, the \emph{theory} of groups as a mathematical object in its own right, models for \lts\ take us back to individual groups as mathematical objects.

\begin{mydefinition}

A \demph{model} for a \lt\ \L\ in a finite-product category $C$ is a finite-product preserving functor
\[\L \tra C.\] 
A map of models is a natural transformation between them.  These form a category $\Mod(\L, C)$.

\end{mydefinition}

\begin{myexample}

Let \L\ be the \lt\ for monoids, and $C = \Set$.  Consider a finite-product preserving functor
\[F\: \L \tra C.\]
Writing $F(1) = A$, we must have $F(k) = A^k$.  Then given any $k$-ary operation, that is, morphism $k \tra 1$ in \L, we get a function
\[A^k \tra A.\]
Functoriality and preservation of products ensures that this is precisely a monoid as expected.  Putting $C = \Top$ gives an underlying \emph{space} $A$ with multiplication given by continuous maps, so we get topological monoids as expected.
\end{myexample}

We now discuss the correspondence between \lts\ and monads, which was hinted at in Example~\ref{ex1}.  This was originally analysed by Linton \cite{lin1}.

\begin{proposition}
 Given a monad $T$ on \Set\ we can construct a \lt\ $\L_T$ as the full subcategory of $\fn{Kl} T\op$ whose objects are those of \F.  Moreover if $T$ is finitary 
\[\cat{Mod}(\L_T, \Set) \catequiv \fn{Alg}T.\]
\end{proposition}

\begin{remark}
 It is worth unravelling this a bit.  Recall in Example~\ref{ex1} we saw that the morphisms $2 \tra 1$ in the \lt\ for monoids were given by all the elements of $T[2]$, where $T$ is the free monoid monad and $[2]$ is a 2-element set.  
\end{remark}

So we see that
\[\begin{array}{rcl}
   \L_T(2,1) &=& \Set(1, T[2]) \\
&=& \fn{Kl}T(1,2) \\
 &=& \fn{Kl}T\op(2,1).
  \end{array}\]
More generally a morphism $k \tra m$ is ``$m$ operations of arity $k$'' ie
\[\begin{array}{rcl}
   \L_T(k,m) &=& \Set([m], T[k]) \\
&=& \fn{Kl}T([m], [k])\\
&=& \fn{Kl}T\op ([k], [m]).
  \end{array}\]
Note that this has finite products because \set\ has coproducts.
Now as we have only used finite sets, we cannot hope to have captured all the behaviour of a general monad on \Set---only the finitary part.  Recall that a finitary functor is one that preserves filtered colimits; on \Set\ this amounts to being entirely determined by its action on finite sets as follows.

\begin{proposition}
 Let $F$ be a functor $\Set \tra \Set$.  Then $F$ is finitary if and only if
\[FX = \int\limits^{[n] \in \finset} \hh{-1.5em} F[n] \times X^n.\]
\end{proposition}

\noi This indicates how we can construct a monad from a \lt.

\begin{proposition} (Linton \cite{lin1})
 Given a \lt\ \L\ we can construct a finitary monad $T_\L$  on \Set\ by
\[T_\L X = \int\limits^{[n] \in \finset} \hh{-1.5em} \L(n,1) \times X^n.\]
\end{proposition}

\noi This gives us a correspondence between \lts\ and finitary monads on \Set.

\begin{theorem}
 The constructions $T \tmapsto \L_T$ and $\L \tmapsto T_\L$ extend to functors exhibiting \Law\ as a full coreflective subcategory of \Mnd, the category of monads on \Set.  Moreover, the essential image of the functor
\[\Law \tra \Mnd\]
is given by the finitary monads, that is, the functor becomes an equivalence
\[\Law \tmap{\catequiv} \Mnd_f\]
where $\Mnd_f$ denotes the full subcategory of finitary monads on \set.
\end{theorem}

This paper can be seen as providing several equivalent definitions of \dl\ for \lt\ that extend the above correspondence.


\section{Distributive laws for monads}\label{dls}

In this work we will be thinking of \dls\ in two ways:
\begin{enumerate}
 \item a way of combining algebraic theories to provide a composite theory, and
\item more generally: an abstract structure giving a way of composing monads to produce a composite monad inside any bicategory $\cB$.
\end{enumerate}
In this section we will simply recall the basic definitions.  None of the material in this section is new.  We first recall the classical theory of distributive laws, stated with respect to monads \emph{on} categories.

\begin{defn} {\upshape (Beck \cite{bec1})} Let $S$ and $T$ be monads on a 
category \cl{C}.  A \demph{distributive law of $S$ over $T$} consists of a 
natural transformation $\lambda\: ST \Rightarrow TS$ such that the following 
diagrams commute.

\[\psset{unit=0.1cm,labelsep=2pt,nodesep=3pt,linewidth=0.8pt}
\pspicture(80,20)


\rput(0,0){
\rput(0,0){\rnode{a}{$ST$}}
\rput(20,0){\rnode{b}{$TS$}}
\rput(10,15){\rnode{c}{$T$}}

\ncline{->}{a}{b}\nbput{\scriptsize $\lambda$}
\ncline{->}{c}{a}\nbput{\scriptsize $\eta^ST$}
\ncline{->}{c}{b}\naput{\scriptsize $T\eta^S$}

}


\rput(40,0){
\rput(0,15){\rnode{a}{$S^2T$}}
\rput(20,15){\rnode{b}{$STS$}}
\rput(40,15){\rnode{c}{$TS^2$}}
\rput(0,0){\rnode{d}{$ST$}}
\rput(40,0){\rnode{e}{$TS$}}

\ncline{->}{a}{b}\naput{\scriptsize $S\lambda$}
\ncline{->}{b}{c}\naput{\scriptsize $\lambda S$}
\ncline{->}{a}{d}\nbput{\scriptsize $\mu^ST$}
\ncline{->}{c}{e}\naput{\scriptsize $T\mu^S$}
\ncline{->}{d}{e}\nbput{\scriptsize $\lambda$}

}

\endpspicture\]

\[\psset{unit=0.1cm,labelsep=2pt,nodesep=3pt,linewidth=0.8pt}
\pspicture(0,-3)(80,18)


\rput(0,0){
\rput(0,0){\rnode{a}{$ST$}}
\rput(20,0){\rnode{b}{$TS$}}
\rput(10,15){\rnode{c}{$S$}}

\ncline{->}{a}{b}\nbput{\scriptsize $\lambda$}
\ncline{->}{c}{a}\nbput{\scriptsize $S\eta^T$}
\ncline{->}{c}{b}\naput{\scriptsize $\eta^T S$}

}


\rput(40,0){
\rput(0,15){\rnode{a}{$ST^2$}}
\rput(20,15){\rnode{b}{$TST$}}
\rput(40,15){\rnode{c}{$T^2 S$}}
\rput(0,0){\rnode{d}{$ST$}}
\rput(40,0){\rnode{e}{$TS$}}

\ncline{->}{a}{b}\naput{\scriptsize $\lambda T$}
\ncline{->}{b}{c}\naput{\scriptsize $T\lambda$}
\ncline{->}{a}{d}\nbput{\scriptsize $S\mu^T$}
\ncline{->}{c}{e}\naput{\scriptsize $\mu^TS$}
\ncline{->}{d}{e}\nbput{\scriptsize $\lambda$}

}

\endpspicture\]

%

\end{defn}

The main theorem about distributive laws tells us about new monads that arise canonically as a result of the distributive law.  In this work we will mostly be interested in the composite monad.

\begin{thm}[Beck, \cite{bec1}]\label{barrwells}

The following are equivalent:
\numarabic
\begin{itemize}
\item A distributive law of $S$ over $T$.
\item A lifting of the monad $T$ to a monad $T'$ on $S\mbox{{\upshape -Alg}}$.
\item An extension of the monad $S$ to a monad $\tilde{S}$ on $\mbox{{\upshape Kl}} (T)$.
\end{itemize}
It follows that $TS$ canonically acquires the structure of a monad, whose category of algebras coincides with that of the lifted monad $T'$, and whose Kleisli category coincides with that of $\tilde{S}$.
\end{thm}

\begin{eg}\label{ring}

{\bfseries (Rings)}

$\cl{C} = \cat{Set}$

$S = \mbox{free monoid monad}$

$T = \mbox{free abelian group monad}$

$\lambda = \mbox{the usual distributive law for multiplication and addition}$ e.g. \[(a+b)(c+d) \mapsto ac+bc+ad+bd.\]

\noi Then the composite monad $TS$ is the free ring monad.
\end{eg}
%
%
%
%
%
%
%

\begin{eg}

{\bfseries (2-categories)}

$\cl{C}=\cat{2-GSet}$, the category of 2-globular sets.

$S = $ monad for vertical composition of 2-cells (1- and 0-cells are unchanged)

$T=$ monad for horizontal composition of 2-cells and 1-cells (0-cells are unchanged)

$\lambda$ is given by the interchange law e.g.
\[
\psset{unit=0.07cm,labelsep=2pt,nodesep=2pt}
\pspicture(0,-5)(120,35)


\rput(21,30){\rnode{c1}{$ST$}}
\rput(97,30){\rnode{c2}{$TS$}}

\rput(21,8){\rnode{d1}{}}
\rput(97,8){\rnode{d2}{}}

\ncline[nodesep=20pt]{->}{c1}{c2}
\ncline[nodesepA=58pt,nodesepB=63pt]{|->}{d1}{d2}

\rput[l](0,0){

\rput(0,0){

\rput(0,10){\rnode{a1}{$\cdot$}}  
\rput(20,10){\rnode{a2}{$\cdot$}}  
\rput(40,10){\rnode{a3}{$\cdot$}}

{
\rput[c](10,12){\psset{unit=1mm,doubleline=true,arrowinset=0.6,arrowlength=0.5
, arrowsize=0.5pt 2.1,nodesep=0pt,labelsep=2pt}
\pcline{->}(0,3)(0,0) \naput{{\scriptsize $$}}}}

{
\rput[c](30,12){\psset{unit=1mm,doubleline=true,arrowinset=0.6,arrowlength=0.5
, arrowsize=0.5pt 2.1,nodesep=0pt,labelsep=2pt}
\pcline{->}(0,3)(0,0) \naput{{\scriptsize $$}}}}

\nccurve[angleA=60,angleB=120,ncurv=1]{->}{a1}{a2}\naput{{\scriptsize $$}}
\ncline{->}{a1}{a2}

\nccurve[angleA=60,angleB=120,ncurv=1]{->}{a2}{a3}\naput{{\scriptsize $$}}
\ncline{->}{a2}{a3}

}

\rput[l](0,-5){

\rput(0,10){\rnode{b1}{$\cdot$}}  
\rput(20,10){\rnode{b2}{$\cdot$}}  
\rput(40,10){\rnode{b3}{$\cdot$}}

\ncline{->}{b1}{b2}
\nccurve[angleA=-60,angleB=-120,ncurv=1]{->}{b1}{b2}\nbput{{\scriptsize $$}}

\ncline{->}{b2}{b3}
\nccurve[angleA=-60,angleB=-120,ncurv=1]{->}{b2}{b3}\nbput{{\scriptsize $$}}

{
\rput[c](10,4.5){\psset{unit=1mm,doubleline=true,arrowinset=0.6,arrowlength=0.5
, arrowsize=0.5pt 2.1,nodesep=0pt,labelsep=2pt}
\pcline{->}(0,3)(0,0) \naput{{\scriptsize $$}}}}

{
\rput[c](30,4.5){\psset{unit=1mm,doubleline=true,arrowinset=0.6,arrowlength=0.5
, arrowsize=0.5pt 2.1,nodesep=0pt,labelsep=2pt}
\pcline{->}(0,3)(0,0) \naput{{\scriptsize $$}}}}

}
}

\rput(75,-2.5){

\rput(0,0){
\rput(0,10){\rnode{a1}{$\cdot$}}  
\rput(20,10){\rnode{a2}{$\cdot$}}  

\nccurve[angleA=60,angleB=120,ncurv=1]{->}{a1}{a2}\naput{{\scriptsize $$}}
\nccurve[angleA=-60,angleB=-120,ncurv=1]{->}{a1}{a2}\nbput{{\scriptsize $$}}

\ncline{->}{a1}{a2}

\rput[c](10,12){
\psset{unit=1mm,doubleline=true,arrowinset=0.6,arrowlength=0.5
, arrowsize=0.5pt 2.1,nodesep=0pt,labelsep=2pt}
\pcline{->}(0,3)(0,0) \naput{{\scriptsize $$}}}

\rput[c](10,4.5){
\psset{unit=1mm,doubleline=true,arrowinset=0.6,arrowlength=0.5
, arrowsize=0.5pt 2.1,nodesep=0pt,labelsep=2pt}
\pcline{->}(0,3)(0,0) \naput{{\scriptsize $$}}}

}

\rput(25,0){
\rput(0,10){\rnode{a1}{$\cdot$}}  
\rput(20,10){\rnode{a2}{$\cdot$}}  

\nccurve[angleA=60,angleB=120,ncurv=1]{->}{a1}{a2}\naput{{\scriptsize $$}}
\nccurve[angleA=-60,angleB=-120,ncurv=1]{->}{a1}{a2}\nbput{{\scriptsize $$}}

\ncline{->}{a1}{a2}

\rput[c](10,12){
\psset{unit=1mm,doubleline=true,arrowinset=0.6,arrowlength=0.5
, arrowsize=0.5pt 2.1,nodesep=0pt,labelsep=2pt}
\pcline{->}(0,3)(0,0) \naput{{\scriptsize $$}}}

\rput[c](10,4.5){
\psset{unit=1mm,doubleline=true,arrowinset=0.6,arrowlength=0.5
, arrowsize=0.5pt 2.1,nodesep=0pt,labelsep=2pt}
\pcline{->}(0,3)(0,0) \naput{{\scriptsize $$}}}

}

}

\endpspicture
\]
\end{eg}

The main theorem of \cite{che17} generalises the notion of distributive law to the case when we have more than two monads interacting with each other, as follows.

\begin{thm}\label{idl}

Fix $n \geq 3$.  Let $T_1, \ldots, T_n$ be monads on a category \cl{C}, equipped with

\begin{itemize}

\item for all $i > j$ a distributive law $\lambda_{ij}: T_i T_j \Rightarrow T_j T_i$, satisfying

\item for all $i>j>k$ the ``Yang-Baxter'' equation given by the commutativity of the following diagram

\begin{equation}\label{yb}
\psset{unit=0.1cm,labelsep=2pt,nodesep=2pt}
\pspicture(0,-3)(80,33)


\rput(0,15){\rnode{a}{$T_iT_jT_k$}}
\rput(15,30){\rnode{b}{$T_jT_iT_k$}}
\rput(45,30){\rnode{c}{$T_jT_kT_i$}}
\rput(60,15){\rnode{d}{$T_kT_jT_i$}}

\rput(15,0){\rnode{e}{$T_iT_kT_j$}}
  \rput(45,0){\rnode{f}{$T_kT_iT_j$}}

\ncline{->}{a}{b}\naput{\scriptsize $\lambda_{ij}T_k$}
\ncline{->}{b}{c}\naput{\scriptsize $T_j\lambda_{ik}$}
\ncline{->}{c}{d}\naput{\scriptsize $\lambda_{jk}T_i$}

\ncline{->}{a}{e}\nbput{\scriptsize $T_i\lambda_{jk}$}
\ncline{->}{e}{f}\nbput{\scriptsize $\lambda_{ik}T_j$}
\ncline{->}{f}{d}\nbput{\scriptsize $T_k\lambda_{ij}$}

\endpspicture
\end{equation}

%
%

\end{itemize}

\noindent Then for all $1\leq i < n$  we have canonical monads
\[T_1 T_2 \cdots T_i \quad \mbox{and} \quad T_{i+1} T_{i+2} \cdots T_n\]

\noindent together with a distributive law of \ $T_{i+1} T_{i+2} \cdots T_n$\  over \ $T_1 T_2 \cdots T_i$\ i.e.
\[(T_{i+1} T_{i+2} \cdots T_n) (T_1 T_2 \cdots T_i) \Rightarrow (T_1 T_2 \cdots T_i)(T_{i+1} T_{i+2} \cdots T_n)\]

\noindent given by the obvious composites of the $\lambda_{ij}$.  Moreover, all the induced monad structures on \ $T_1 T_2 \cdots T_n$ are the same.

\end{thm}

\begin{defn} \label{dsm}
A  \demph{distributive series of $n$ monads} is a system of monads and distributive laws as in Theorem~\ref{idl}.
\end{defn}

\begin{eg} Rings can be constructed from the following 
distributive series of 3 monads on \set.
\[\begin{array}{ccl}
A &=& \mbox{monad for associative non-unital binary multiplication $\times$} \\
B &=& \mbox{monad for pointed sets i.e. $X \mapsto X \coprod \{1\}$} \\
C &=& \mbox{free additive abelian group monad}
\end{array}\]
\end{eg}

\begin{myexample} 
Strict $n$-categories can be constructed from a distributive series of $n$ monads on $n$-globular sets, as a generalisation of the 2-category case.  Here there is a monad $T_i$ for each $0 \leq i \leq n-1$ giving ``composition along bounding $n$-cells.
\end{myexample}

In his classic paper {\em The formal theory of monads} \cite{str1} Street defines for any 2-category \cl{B} a 2-category $\cat{Mnd}(\cl{B})$ of monads in \cl{B}.  Then distributive laws arise as monads in $\cat{Mnd}(\cl{B})$.  While we will not use that particular, and appealing, fact, we will certainly be looking at monads and distributive laws inside various 2-categories and in fact bicategories, which can be done by invoking appropriate coherence conditions and results.


\section{Monads in profunctors}\label{prof}

In this section we give the most straightforward but perhaps least intuitive 
definition of distributive laws for \lts.  We start to make use of the notion of a monad \emph{in} a bicategory. We use the bicategory \Prof\ of profunctors, and simply observe that all \lts\ are monads on $\F\op$ in \Prof\ (though 
not all monads on $\F\op$ are \lts); this result and those leading up to it are 
standard.  We can thus simply look at \dls\ between these monads.  It is not 
immediately obvious why this should be the right definition and we will defer this 
justification to the last section.  

First we set our notational conventions.

\begin{mydefinition}
 We write \Prof\ for the bicategory given as follows.

\begin{itemize}
 \item 0-cells are small categories,
\item a 1-cell $\C \pmap{F} \D$ is a functor $\D\op \times \C \stra \Set$,
\item 2-cells are natural transformations. 
\end{itemize}

Composition of profunctors $\C \pmap{F} \D \pmap{G} \E$ is by the usual coend formula
\[(G \circ F)(e,c) = \int\limits^{d \in \D} \hh{-0.4em} G(e,d) \times F(d,c)\]
and is only weakly associative and unital.

\end{mydefinition}

Profunctors turn out to be the same as bimodules internal to the bicategory of spans.  This fact will be useful to us both technically and conceptually in Section~\ref{factsystems}.

\begin{mydefinition}
 We write \Span\ for the bicategory of spans given as follows.

\begin{itemize}
 \item 0-cells are sets,
\item a 1-cell $C \dmap{X} D$ is a span 
%
%
%
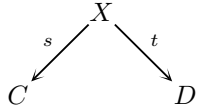
\[
\psset{unit=0.1cm,labelsep=3pt,nodesep=1.5pt}
\pspicture(30,15)


\rput(11,11){\rnode{a1}{$X$}}  
\rput(0,0){\rnode{a3}{$C$}}  
\rput(22,0){\rnode{a2}{$D$}}  

\ncline{->}{a1}{a3} \nbput{{\scriptsize $s$}} 
\ncline{->}{a1}{a2} \naput{{\scriptsize $t$}} 
\endpspicture\]

\item 2-cells are morphisms of spans.
\end{itemize}
Composition of 1-cells is by pullback, so given $C \dmap{X} D \dmap{Y} E$ we have 
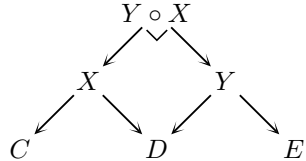
\[
\psset{unit=0.09cm,labelsep=3pt,nodesep=2pt}
\pspicture(-20,0)(40,22)



\rput(10,20){\rnode{a1}{$Y\circ X$}}  
\rput(0,10){\rnode{a3}{$X$}}  
\rput(20,10){\rnode{a2}{$Y$}}  
\rput(10,0){\rnode{a4}{$D$}}  
\rput(-10,0){\rnode{a5}{$C$}}  
\rput(30,0){\rnode{a6}{$E$}} 

\ncline{->}{a1}{a3} \nbput{{\scriptsize $$}} 
\ncline{->}{a1}{a2} \naput{{\scriptsize $$}} 
\ncline{->}{a3}{a4} \nbput{{\scriptsize $$}} 
\ncline{->}{a2}{a4} \naput{{\scriptsize $$}} 

\ncline{->}{a3}{a5} \nbput{{\scriptsize $$}} 
\ncline{->}{a2}{a6} \nbput{{\scriptsize $$}} 

\psline[linewidth=0.6pt]{-}(8.5,17.5)(10,16)(11.5,17.5)

\endpspicture
\]

\end{mydefinition}

\begin{mydefinition}
Given any bicategory $K$ and monads $X, Y$ inside it, a $(Y,X)$-bimodule $A$ is given by a 1-cell $x \tmap{A} y$ in $K$ equipped with 2-cell actions 
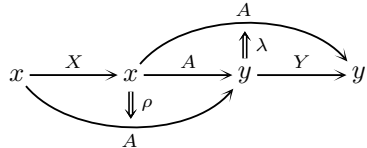
\[
\psset{unit=0.1cm,labelsep=2pt,nodesep=2pt}
\pspicture(45,20)

\rput(0,10){\rnode{a1}{$x$}}  
\rput(15,10){\rnode{a2}{$x$}}  
\rput(30,10){\rnode{a3}{$y$}}  
\rput(45,10){\rnode{a4}{$y$}}  

\ncline{->}{a1}{a2} \naput{{\scriptsize $X$}}
\ncline{->}{a2}{a3} \naput{{\scriptsize $A$}}
\ncline{->}{a3}{a4} \naput{{\scriptsize $Y$}}

\pcline[linewidth=0.6pt,doubleline=true,arrowinset=0.6,arrowlength=0.8,arrowsize=0.5pt 2.1]{->}(15,8.5)(15,3.5)\naput{{\scriptsize $\rho$}}
\pcline[linewidth=0.6pt,doubleline=true,arrowinset=0.6,arrowlength=0.8,arrowsize=0.5pt 2.1]{->}(30,11.5)(30,16.5)\nbput{{\scriptsize $\lambda$}}

\ncarc[arcangle=50]{->}{a2}{a4}\naput{{\scriptsize $A$}}
\ncarc[arcangle=-50]{->}{a1}{a3}\nbput{{\scriptsize $A$}}

\endpspicture
\]
satisfying the usual bimodule axioms: $\lambda$ is compatible with the structure of $X$, $\rho$ with the structure of $Y$ and $\lambda$ and $\rho$ with each other.
\end{mydefinition}

\noi Provided $K$ has enough structure, bimodules are the 1-cells of a bicategory as follows.

\begin{mydefinition}
Let $K$ be a bicategory with coequalisers of 2-cells that are preserved by 
left and right composition with 1-cells. We write $\Mod(K)$ for the bicategory 
of bimodules in $K$, given as follows.

\begin{itemize}
 \item 0-cells are the monads in $K$,
\item a 1-cell $X \cmap{A} Y$ is a $(Y,X)$-bimodule (note direction).

%
%
%
%
%
%
%
%
%
%
%

\item 2-cells are bimodule maps.

\item Composition of 1-cells is by coequaliser: given
\[X \cmap{A} Y \cmap{B} Z\]
given by 1-cells
\[x \tmap{X} x \tmap{A} y \tmap{Y} y \tmap{B} z \tmap{Z} z\]
we take the coequaliser
\[
\psset{labelsep=2pt}
\pspicture(70,7)
\rput(0,3){\rnode{A}{$B \circ Y \circ A$}}
\rput(27,3){\rnode{B}{$B \circ A$}}
\rput(50,3){\rnode{C}{$B \otimes_Y A$}}
\psset{nodesep=3pt,arrows=->}
\ncline[offset=3pt]{A}{B}\naput{$\rho \circ A$}
\ncline[offset=-3pt]{A}{B}\nbput{$B \circ \lambda$}
\ncline{B}{C}\naput{$$}
\endpspicture
\]

\noi $B \otimes_Y A$ is then the composite $(Z,X)$-bimodule required.

\end{itemize}

\end{mydefinition}

\noi Combining these two constructions gives another way of thinking 
of profunctors, with some care over dualities.  

\begin{myexample}
 The bicategory $\Mod(\Span)$ is given as follows.

\begin{itemize}
 \item 0-cells are monads in \Span\, that is, small categories.

\item Given categories $X, Y$ with underlying spans
%
%
\[
\psset{unit=0.1cm,labelsep=3pt,nodesep=1.5pt}
\pspicture(60,10)

\rput[b](0,0){
\pspicture(30,15)


\rput(9,9){\rnode{a1}{$X_1$}}  
\rput(0,0){\rnode{a3}{$X_0$}}  
\rput(18,0){\rnode{a2}{$X_0$}}  

\ncline{->}{a1}{a3} \nbput{{\scriptsize $s$}} 
\ncline{->}{a1}{a2} \naput{{\scriptsize $t$}} 
\endpspicture}

%
%
\rput[b](35,0){
\pspicture(30,15)


\rput(9,9){\rnode{a1}{$Y_1$}}  
\rput(0,0){\rnode{a3}{$Y_0$}}  
\rput(18,0){\rnode{a2}{$Y_0$}}  

\ncline{->}{a1}{a3} \nbput{{\scriptsize $s$}} 
\ncline{->}{a1}{a2} \naput{{\scriptsize $t$}} 
\endpspicture}

\endpspicture\]

\noi a 1-cell $X \dmap{A} Y$ has underlying span of the form

%
%
\[
\psset{unit=0.1cm,labelsep=3pt,nodesep=1.5pt}
\pspicture(30,10)


\rput(9,9){\rnode{a1}{$A_1$}}  
\rput(0,0){\rnode{a3}{$X_0$}}  
\rput(18,0){\rnode{a2}{$Y_0$}}  

\ncline{->}{a1}{a3} \nbput{{\scriptsize $s$}} 
\ncline{->}{a1}{a2} \naput{{\scriptsize $t$}} 
\endpspicture\]

The elements of $A$ can be thought of as arrows with source in $X$ and target in $Y$.  The left $Y$-action is a map of spans

\[
\psset{unit=0.09cm,labelsep=1pt,nodesep=2pt}
\pspicture(50,20)

\rput[b](0,0){\psset{unit=0.08cm}\pspicture(-20,0)(40,22)



\rput(10,20){\rnode{a1}{$.$}}  
\rput(0,10){\rnode{a3}{$A_1$}}  
\rput(20,10){\rnode{a2}{$Y_1$}}  
\rput(-10,0){\rnode{a5}{$X_0$}}  
\rput(10,0){\rnode{a4}{$Y_0$}}  
\rput(30,0){\rnode{a6}{$Y_0$}} 

\ncline{->}{a1}{a3} \nbput{{\scriptsize $$}} 
\ncline{->}{a1}{a2} \naput{{\scriptsize $$}} 

\ncline{->}{a3}{a5} \nbput{{\scriptsize $s$}} 
\ncline{->}{a3}{a4} \naput{{\scriptsize $t$}} 
\ncline{->}{a2}{a4} \nbput{{\scriptsize $s$}} 
\ncline{->}{a2}{a6} \naput{{\scriptsize $t$}} 

\psline[linewidth=0.6pt]{-}(8.5,17.5)(10,16)(11.5,17.5)

\endpspicture}

\pcline[linewidth=1.5pt]{->}(25,10)(40,10)

\rput[b](60,3){\psset{unit=0.12cm}\pspicture(20,12)


\rput(9,9){\rnode{a1}{$A_1$}}  
\rput(0,0){\rnode{a3}{$X_0$}}  
\rput(18,0){\rnode{a2}{$Y_0$}}  

\ncline{->}{a1}{a3} \nbput{{\scriptsize $s$}} 
\ncline{->}{a1}{a2} \naput{{\scriptsize $t$}} 
\endpspicture}
\endpspicture
\]
giving us a way of post-composing arrows in $A$ with those of $Y$; the module axioms tell us that this respects composition in $Y$.  Similarly for the left $X$-action.  The left-right compatibility then gives us associativity for composing three arrows
\[\tmap{\in X} \tmap{\in A} \tmap{\in Y}.\]
\end{itemize}

\end{myexample}

\noi In \cite{ben2} B\'enabou first defines profunctors 
(``distributeurs'') directly 
as functors $\D\op \times \C \stra \Set$. He then defines profunctors internal to a 
bicategory $\cl{E}$ as bimodules in the bicategory $\fn{Span} \cE$ of spans 
internal to 
$\cl{E}$ as follows.

\begin{mydefinition}\label{ben} \cite{ben2}
Let $\cE$ be a category with pullbacks and coequalisers that commute.  
\begin{itemize}
 \item Write $\fn{Span} \cE$ for the bicategory of spans in $\cE$.
\item Define $\fn{Prof} \cE$ to be the bicategory $\Mod(\fn{Span} \cE)\op$.  
Thus 0-cells are monads in $\fn{Span} \cE$, that is, categories internal to 
$\cE$.
\end{itemize}

\end{mydefinition}

\begin{myremarks} \ \\[-1em]
\begin{enumerate}
 \item We need pullbacks to define composition of spans, and we need the 
coequaliser condition to define composition of profunctors.

\item We need to take the dual here for reasons that will become clear later..


\end{enumerate}

\end{myremarks}

Thus according to this approach profunctors in \Set\ are bimodules in \Span\ by 
definition.   Although not stated it seems clear that the 
intention is for profunctors in $\cl{E}$ to be a generalisation of basic profunctors 
in the sense that the notions coincide in the case $\cl{E} =  \Set.$  This is the 
content of the following proposition.

\begin{myproposition}
 There is a biequivalence of bicategories
\[\Prof\op \catequiv \Mod(\Span).\]
\end{myproposition}

\begin{prf} (Sketch.) First we construct a functor
\[\Prof\op \tra \Mod(\Span).\]
The 0-cells on both sides are small categories, thus we set the action of the 
functor on 0-cells to be the identity.  

For the action on 1-cells, we start with a profunctor $Y\pmap{F} X$, that is a
functor $X\op \times Y \tmap{F} \Set$, and construct a bimodule $X \cra Y$, that is, 
a $(Y,X)$-bimodule, as follows.  First take the underlying span $X \dra Y$ to be:
\[
\psset{unit=0.1cm,labelsep=3pt,nodesep=1.5pt}
\pspicture(40,15)


\rput(13,13){\rnode{a1}{\scalebox{0.8}{$\displaystyle\coprod_{x,y}F(x,y)$}}}  
\rput(0,0){\rnode{a3}{{\scriptsize $X_0$}}}  
\rput(26,0){\rnode{a2}{{\scriptsize $Y_0$}}}  

\ncline{->}{a1}{a3} \nbput{{\scriptsize $s$}} 
\ncline{->}{a1}{a2} \naput{{\scriptsize $t$}} 
\endpspicture\]

\noi The left $Y$- and right $X$-actions are given by the actions of $F$ on 
morphisms as follows.  For the $Y$-action we need a map of spans

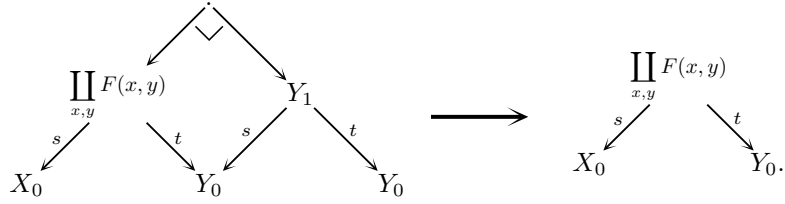
\[
\psset{unit=0.09cm,labelsep=1pt,nodesep=1pt}
\pspicture(50,25)

\rput[b](-10,0){\psset{unit=0.12cm}\pspicture(-20,0)(40,22)



\rput(10,20){\rnode{a1}{$.$}}  
\rput(0,10){\rnode{a3}{\scalebox{0.8}{$\displaystyle\coprod_{x,y}F(x,y)$}}}  
\rput(20,10){\rnode{a2}{$Y_1$}}  
\rput(-10,0){\rnode{a5}{$X_0$}}  
\rput(10,0){\rnode{a4}{$Y_0$}}  
\rput(30,0){\rnode{a6}{$Y_0$}} 

\ncline{->}{a1}{a3} \nbput{{\scriptsize $$}} 
\ncline{->}{a1}{a2} \naput{{\scriptsize $$}} 

\ncline{->}{a3}{a5} \nbput{{\scriptsize $s$}} 
\ncline{->}{a3}{a4} \naput{{\scriptsize $t$}} 
\ncline{->}{a2}{a4} \nbput{{\scriptsize $s$}} 
\ncline{->}{a2}{a6} \naput{{\scriptsize $t$}} 

\psline[linewidth=0.6pt]{-}(8.5,17.5)(10,16)(11.5,17.5)

\endpspicture}

\pcline[linewidth=1.5pt]{->}(22,10)(37,10)

\rput[b](60,3){\psset{unit=0.13cm}\pspicture(20,12)


\rput(9,9){\rnode{a1}{\scalebox{0.8}{$\displaystyle\coprod_{x,y}F(x,y)$}}}  
\rput(0,0){\rnode{a3}{$X_0$}}  
\rput(18,0){\rnode{a2}{$Y_0.$}}  

\ncline{->}{a1}{a3} \nbput{{\scriptsize $s$}} 
\ncline{->}{a1}{a2} \naput{{\scriptsize $t$}} 
\endpspicture}
\endpspicture
\]

\vv{1em}

\noi An element in the pullback is a pair $(\alpha \in F(x,y), f\in Y_1(y,y'))$.  Now we have
\[Ff\: F(x,y) \tra F(x,y')\]
so we define the action by
\[(\alpha, f) \tmapsto Ff(\alpha).\]
The $X$-action is constructed similarly.

Now we construct a functor
\[\Mod(\Span) \tra \Prof\op\]
which again is the identity on 0-cells.  Given categories $X, Y$ and a bimodule $X \cmap{A} Y$, that is, a $(Y,X)$-bimodule with underlying span
%
%
\[
\psset{unit=0.1cm,labelsep=3pt,nodesep=1.5pt}
\pspicture(20,12)


\rput(9,9){\rnode{a1}{$A_1$}}  
\rput(0,0){\rnode{a3}{$X_0$}}  
\rput(18,0){\rnode{a2}{$Y_0,$}}  

\ncline{->}{a1}{a3} \nbput{{\scriptsize $s$}} 
\ncline{->}{a1}{a2} \naput{{\scriptsize $t$}} 
\endpspicture\]
say, we construct a profunctor $Y \pra X$, that is, a functor
\[X\op \times Y \tmap{F} \Set,\]
by $F(x,y) = A(x,y)$, that is, the pre-image in $A_1$ of the pair $(x,y)$.  
Functoriality comes from the left and right actions. It is routine to check 
that this gives a biequivalence of bicategories. \end{prf}

\begin{myremark} \   Note that when we discuss factorisation systems in 
Section~\ref{factsystems} it is useful to think in terms of spans, but for the 
comparison in Section~\ref{comparison} it is useful to think in terms of 
profunctors.
\end{myremark}

We are going to show that Lawvere theories arise as certain monads in \Prof.  In fact 
the monads in \Prof\ are \emph{any} identity-on-objects functors.  This is fairly 
easy to prove directly, but it is also a special case of the following standard 
result.


\begin{theorem}\label{theoremA}\label{thmA}
Let $K$ be a bicategory, $x$ a 0-cell, and $X$ a monad on $x$.  Then there is an equivalence of categories
\[ \Mon\big( ( \Mod K ) (X,X) \big) \catequiv X/ \Mon\big( K(x,x) \big).\] 

\end{theorem}

\noi Note that here we write $\Mon \cV$ for the category of monoids in a monoidal 
category $\cV$, and $\cB(b,b)$ for the monoidal category of 1-cells $b \tra b$ in 
a bicategory $\cB$.  Thus on the left hand side we 

\begin{enumerate}
 \item form the bicategory of bimodules in $K$,
\item take the monoidal category of 1-cells $X \cra X$ in this bicategory, and
\item take the category of monoids in this monoidal category.
\end{enumerate}

\noi For the right hand side we 
\begin{enumerate}
 \item take the monoidal category of 1-cells $x \tra x$ in $K$,
\item take the category of monoids in this monoidal category, and
\item slice this category under $X$.
\end{enumerate}

\begin{corollary}
A monad in $\Mod(\Span)$ on $X$ consists of a category $A$ and an 
identity-on-objects functor $X \tra A$.
\end{corollary}

\begin{prf} $X$ is a 0-cell of $\Mod(\Span)$ so is a monad in \Span\, that is, a 
small category with object set $x$, say.   Now a monad in 
$\Mod(\Span)$ on $X$ is a monoid in $\Mod((\Span)(X,X))$ by definition, so by 
Theorem~\ref{theoremA} it is an object of $X/\Mon(\Span(x,x))$.  Now

\begin{itemize}
 \item a monoid in $\Span(x,x)$ is a category with the same objects as $X$, and
\item a morphism of monoids in $\Span(x,x)$ is an identity-on-objects functor.
\end{itemize}
So the objects of $X/\Mon(\Span(x,x))$ are precisely identity-on-objects 
functors $X \tra A$.    
\end{prf}

\begin{corollary}
A monad in $\Prof\op$ on $X$ consists of a category $A$ and an 
identity-on-objects functor $X \tra A$.
\end{corollary}

\begin{myremark}
 It is illuminating to sketch a direct proof of this result.  A monad $X \prof X$ in $\Mod(\Span)$ is an $(X,X)$-bimodule that is also a monad.  That is, it has a left and right $X$-action but also a unit and multiplication of its own.  Note that $X$ is itself a monad in \Span, with underlying span
%

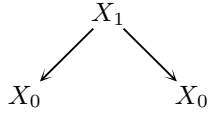
\[
\psset{unit=0.1cm,labelsep=3pt,nodesep=1.5pt}
\pspicture(30,12)


\rput(11,11){\rnode{a1}{$X_1$}}  
\rput(0,0){\rnode{a3}{$X_0$}}  
\rput(22,0){\rnode{a2}{$X_0$}}  

\ncline{->}{a1}{a3} \nbput{{\scriptsize $$}} 
\ncline{->}{a1}{a2} \naput{{\scriptsize $$}} 
\endpspicture\]
say.  So for the monad $X \prof X$ we have a span on the same objects as $X$, say 
%

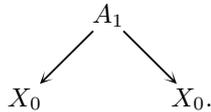
\[
\psset{unit=0.1cm,labelsep=3pt,nodesep=1.5pt}
\pspicture(30,12)


\rput(11,11){\rnode{a1}{$A_1$}}  
\rput(0,0){\rnode{a3}{$X_0$}}  
\rput(22,0){\rnode{a2}{$X_0$.}}  

\ncline{->}{a1}{a3} \nbput{{\scriptsize $$}} 
\ncline{->}{a1}{a2} \naput{{\scriptsize $$}} 
\endpspicture\]
Essentially
\begin{itemize}
 \item the monad structure makes this into a category $A$, say,
\item the left/right $X$-actions tell us how to map $X_1$ to $A_1$, 
\item the way composition of bimodules works ensures that the composition of $A$ is 
compatible with that of $X$, that is, that we have a \emph{functor} $X \tra A$.
\end{itemize}
Similarly we can sketch a direct proof of the result in $\Prof\op$:  given a monad 
$A:\C\op \times \C \stra \Set$ we get a 
category $\bb{A}$ by setting $\bb{A}(a,b) = A(a,b)$ and using the unit and 
multiplication of the monad to give identities and composition.  To construct an 
identity-on-objects functor $\C \stra \bb{A}$ we use the functoriality of $A$ 
which has the effect of producing left and right actions of the morphisms of $\C$ on 
the morphisms of $\bb{A}$; taking the action on identities then gives the functor.

\end{myremark}

\begin{cor} \label{ltmonadinprof}
 Every Lawvere theory $\F\op \tmap{\alpha_A} A$ is a monad on $\Fop$ in $\Prof\op$. 
Conversely a monad $\Fop$ in $\Prof\op$ is a category $A$ equipped 
with an identity-on-objects functor $\F\op \tmap{\alpha_A} A$; it is a Lawvere 
theory precisely if the category $A$ has finite products and the functor $\alpha_A$ 
preserves them.
\end{cor}

\begin{myremark}
 At this point it might seem that we should have started with the opposite (dual) definition of \Prof, which is also standard (and equivalent).  However, in Section~\ref{kleisli} we cannot use that version.
\end{myremark}

\noi Although not every monad on $\F\op$ in $\Prof\op$ is a Lawvere theory, given two Lawvere theories expressed in this way, we can define distributive laws between them.

\begin{mydefinition}\label{def1} ({\scp{``prof''})}
 Given Lawvere theories $A$ and $B$, a \demph{distributive law} of $A$ over $B$ is a distributive law of $A$ over $B$ expressed as monads in $\Prof\op$.  Iterated distributive laws are defined likewise, as in Theorem~\ref{idl}. \end{mydefinition}

\begin{proposition}\label{compositeprof}
 The resulting composite monad $BA$ is also a Lawvere theory.
\end{proposition}

\noi Note that the issue here is finite products---\emph{a priori} our distributive law makes $BA$ into a monad on $\F\op$ in $\Prof$, that is, an identity-on-objects functor $\F\op \tra BA$; for this to be a Lawvere theory we need to prove that $BA$ has finite products and that the functor preserves them.  We defer this proof, and further justification of the definition, until Section~\ref{comparison} (Corollary~\ref{compprof}), as the comparison proceeds via the definitions that we will introduce in subsequent sections.

In the next section we give a more explicit characterisation of such a \dl, using the language of factorisation systems.


\section{Factorisation systems}\label{factsystems}

We will use a notion of factorisation system as given by Rosebrugh and Wood in 
\cite{rw1}, but slightly more general.  Some stages of generalisation of 
notions of factorisation system can be seen as follows:

\begin{enumerate}
 \item Strict factorisation systems on a category $C$.
\item Orthogonal factorisation systems on $C$.
\item Factorisation systems over $I$ where $I$ is a subgroupoid of $C$ 
\cite{rw1}; orthogonal factorisation systems are a special case.
\item Factorisation systems over $J$ where $J$ is a subcategory of $C$.
\end{enumerate}

\noi We include some basic definitions here as the terminology in the literature 
is not entirely uniform. There are also many equivalent formulations; a helpful 
exposition can be found in \cite{rie3}

\begin{mydefinition}

A \demph{strict \fs} on a category $C$ is a pair $(L,R)$ of subcategories of $C$, with the same objects as $C$ (lluf), such that every morphism of $C$ can be factorised \emph{uniquely} as a composite
\[\tmap{l} \tmap{r}\]
with $l \in L$ and $r \in R$.

\end{mydefinition}

\begin{myremarks} \ \\[-1em]
\begin{enumerate}
 \item The uniqueness implies that the intersection of $L$ and $R$ must contain only the identities.
\item It follows that $L \perp R$. (That is, every map in $L$ has the \emph{unique} left lifting property against every map in $R$, and every map in $R$ has the \emph{unique} right lifting property against every map in $L$; this means that lifts exist and are unique.)
\end{enumerate}
\end{myremarks}

\begin{mydefinition}
 An \demph{orthogonal \fs} or simply \demph{\fs} on a category $C$ is a pair $(L,R)$ of lluf subcategories of $C$ containing all isomorphisms, such that every morphism of $C$ can be factorised as a composite
\[\tmap{l} \tmap{r}\]
with $l \in L$ and $r \in R$, uniquely up to unique isomorphism.
\end{mydefinition}

\begin{myremarks} \ \\[-1em]
\begin{enumerate}
\item $L \cap R$ must contain all isomorphisms, so if $C$ contains non-trivial isomorphisms, a strict \fs\ on it is not an orthogonal \fs.
\item It follows that $L \perp R$ and in fact $L = {}^{\perp}R =  {}^{\pitchfork} R$ and $R = L^\perp = L^{\pitchfork}$. Here we write $L^\pitchfork$ for the collection of maps with the right lifting property against all those in $L$, and $L^\perp$ for the collection of maps with the \emph{unique} right lifting property against all those in $L$.  Similarly for $L = {}^{\perp}R$ and  ${}^{\pitchfork} R$ for left liftings.

\end{enumerate}

\end{myremarks}

\begin{myexamples} \ \\[-1em]
 \begin{enumerate}
  \item The pair $(\{\mbox{epi}\}, \{\mbox{mono}\})$ is an orthogonal \fs\ on \Set.
\item The pair $(\{\mbox{bijective-on-objects}\}, \{\mbox{full and faithful}\})$ is an orthogonal \fs\ on \Cat.
\item The pair $(\{\mbox{bijective-on-objects and full}\}, \{\mbox{faithful}\})$ is another orthogonal \fs\ on \Cat.
 \end{enumerate}

\end{myexamples}

\noi There are many naturally-arising \fss\ that are not strict, but the following characterisation by Rosebrugh and Wood \cite{rw1} makes the strict ones of abstract interest.

\begin{theorem}
 Strict \fss\ are precisely \dls\ in \Span.   That is, given a (small) category $C$, a strict factorisation system $(A,B)$ on it is precisely a pair of monads $A$ and $B$ in \Span\ together with a distributive law of $A$ over $B$ such that the composite monad $BA$ is the category $C$.\end{theorem}

\noi Another way of putting this is that a strict \fs\ on a category $C$ is a decomposition of $C$ as a monad in \Span\ into a composite $BA$ via a \dl.

\begin{myremark}  It is worth unravelling this a bit. The composite $BA$ is a 
pullback.  Writing the underlying spans of $A$ and $B$ as
%
%
\[
\psset{unit=0.1cm,labelsep=3pt,nodesep=1.5pt}
\pspicture(30,15)


\rput(11,11){\rnode{a1}{$A_1$}}  
\rput(0,0){\rnode{a3}{$X$}}  
\rput(22,0){\rnode{a2}{$X$}}  

\ncline{->}{a1}{a3} \nbput{{\scriptsize $$}} 
\ncline{->}{a1}{a2} \naput{{\scriptsize $$}} 
\endpspicture
\mbox{and}
%
%
\pspicture(-7,0)(30,15)


\rput(11,11){\rnode{a1}{$B_1$}}  
\rput(0,0){\rnode{a3}{$X$}}  
\rput(22,0){\rnode{a2}{$X$}}  

\ncline{->}{a1}{a3} \nbput{{\scriptsize $$}} 
\ncline{->}{a1}{a2} \naput{{\scriptsize $$}} 
\endpspicture\]
the composite $BA$ is the pullback
\[
\psset{unit=0.09cm,labelsep=3pt,nodesep=2pt}
\pspicture(-20,0)(40,15)



\rput(10,20){\rnode{a1}{$.$}}  
\rput(0,10){\rnode{a3}{$A_1$}}  
\rput(20,10){\rnode{a2}{$B_1$}}  
\rput(10,0){\rnode{a4}{$X$}}  
\rput(-10,0){\rnode{a5}{$X$}}  
\rput(30,0){\rnode{a6}{$X$}} 

\ncline{->}{a1}{a3} \nbput{{\scriptsize $$}} 
\ncline{->}{a1}{a2} \naput{{\scriptsize $$}} 
\ncline{->}{a3}{a4} \nbput{{\scriptsize $$}} 
\ncline{->}{a2}{a4} \naput{{\scriptsize $$}} 

\ncline{->}{a3}{a5} \nbput{{\scriptsize $$}} 
\ncline{->}{a2}{a6} \nbput{{\scriptsize $$}} 

\psline[linewidth=0.6pt]{-}(8.5,17.5)(10,16)(11.5,17.5)

\endpspicture
\]
and is not \emph{a priori} a category.  It consists of pairs of composable morphisms
\[\tmap{\in A} \tmap{\in B}.\]
The \dl\ $AB \tra BA$ tells us how to re-express a composite
\[\tmap{\in B} \tmap{\in A}\]
as one in the ``canonical form''
\[\tmap{\in A} \tmap{\in B}.\]
This makes $BA$ into a category as we can now compose its morphisms: a composable pair in $BA$ will be a composable quadruple
\[\map{\in A} \tmap{\in B} \tmap{\in A} \tmap{\in B}\]
and its composite is obtained by using the distributive law to 
re-express the middle pair to get a string
\[\map{\in A} \tmap{\in A} \tmap{\in B} \tmap{\in B}\]
and then composing in $A$ and in $B$ separately to get a morphism in $BA$.

Note that morphisms in $BA$ are \emph{uniquely} expressible in the form
\[\tmap{\in A} \tmap{\in B}\]
by construction, as these are precisely the morphisms in the pullback.

\end{myremark}

\begin{myexample}{\bfseries(Non-example)} It is instructive to note that this is not the notion we want for \dls\ of \lts.  Let 
\[\alpha \: \F\op \tra A\]
be the \lt\ for (multiplicative) monoids and
\[\beta \: \F\op \tra B\]
be the \lt\ for (additive) Abelian groups.  Thus $X = \ob \F$ in the span notation of the previous remark.  We will now see that $BA$ does not give us the composite theory we want, namely, the theory of rings.

Consider the 3-ary operation $ab + c$ in the theory of rings.  This certainly can be expressed as a composite
\[\tmap{\in A} \tmap{\in B}\]
via
\[ 3 \ltmap{\{ab, c\}} 2 \ltmap{x+y} 1.\]
However, this factorisation is not unique; for example we could also have
\[ 3 \ltmap{\{ab, c, abc\}} 3 \ltmap{x+y} 1.\]
where the first operation adds in a redundant operation $abc$ and the second one forgets it.  Now the two are related via a projection in \Fop\ making the following diagram commute, in the sense that the left-hand triangle commutes in $A$ and the right-hand triangle commutes in $B$. 
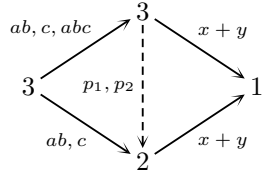
\[
\psset{unit=0.1cm,labelsep=2pt,nodesep=2pt}
\pspicture(0,0)(40,20)



\rput(10,20){\rnode{a1}{$3$}}  
\rput(-5,10){\rnode{a3}{$3$}}  
\rput(25,10){\rnode{a2}{$1$}}  
\rput(10,0){\rnode{a4}{$2$}}  

\ncline{->}{a3}{a1} \naput{{\scriptsize $ab, c, abc$}} 
\ncline{->}{a3}{a4} \nbput{{\scriptsize $ab, c$}} 
\ncline{->}{a1}{a2} \naput{{\scriptsize $x+y$}} 
\ncline{->}{a4}{a2} \nbput{{\scriptsize $x+y$}} 

\ncline[linestyle=dashed, dash=3pt 2pt]{->}{a1}{a4} \nbput{{\scriptsize $p_1, p_2$}} 
\endpspicture\]
However the projection is not an isomorphism, so the factorisation is not unique up to isomorphism. The lesson is that we only want factorisations to be unique up to morphisms in \Fop\ somehow---in fact they are only unique up to zigzags in \Fop.  For example the following two factorisations of the operation $a^2 + a^2$ cannot be related by a single morphism in \Fop: 
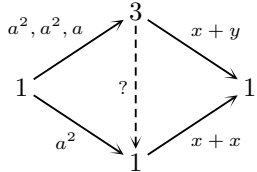
\[
\psset{unit=0.1cm,labelsep=2pt,nodesep=2pt}
\pspicture(0,0)(40,20)



\rput(10,20){\rnode{a1}{$3$}}  
\rput(-5,10){\rnode{a3}{$1$}}  
\rput(25,10){\rnode{a2}{$1$}}  
\rput(10,0){\rnode{a4}{$1$}}  

\ncline{->}{a3}{a1} \naput{{\scriptsize $a^2, a^2, a$}} 
\ncline{->}{a3}{a4} \nbput{{\scriptsize $a^2$}} 
\ncline{->}{a1}{a2} \naput{{\scriptsize $x+y$}} 
\ncline{->}{a4}{a2} \nbput{{\scriptsize $x+x$}} 

\ncline[linestyle=dashed, dash=3pt 2pt]{->}{a1}{a4} \nbput{{\scriptsize ?}} 
\endpspicture\]

\noi We will now show that there is no single morphism in \Fop\ in either direction ($3 \tra 1$ or $1 \tra 3$) that makes the diagram commute.

\begin{itemize}
 \item For morphisms $3 \tra 1$, the only such morphisms are the three projections.  These will clearly not make the resulting right-hand triangle commute.

\item For morphisms $1 \tra 3$, the only such map is the diagonal $\{x,x,x\}$.  This will not make the resulting left-hand triangle commute.

\end{itemize}

\noi So in fact we need a zigzag: 

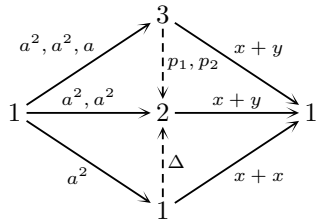
\[
\psset{unit=0.13cm,labelsep=1pt,nodesep=2pt}
\pspicture(0,0)(40,20)



\rput(10,20){\rnode{a1}{$3$}}  
\rput(-5,10){\rnode{a3}{$1$}}  
\rput(25,10){\rnode{a2}{$1$}}  
\rput(10,0){\rnode{a4}{$1$}}  

\ncline{->}{a3}{a1} \naput{{\scriptsize $a^2, a^2, a$}} 
\ncline{->}{a3}{a4} \nbput{{\scriptsize $a^2$}} 
\ncline{->}{a1}{a2} \naput{{\scriptsize $x+y$}} 
\ncline{->}{a4}{a2} \nbput{{\scriptsize $x+x$}} 

\rput(10,10){\rnode{c}{$2$}} 

\ncline{->}{a3}{c} \naput{{\scriptsize $a^2, a^2$}} 
\ncline{->}{c}{a2} \naput{{\scriptsize $x+y$}} 

\ncline[linestyle=dashed, dash=3pt 2pt]{->}{a1}{c} \naput{{\scriptsize $p_1, p_2$}} 

\ncline[linestyle=dashed, dash=3pt 2pt]{->}{a4}{c} \nbput{{\scriptsize $\Delta$}} 
\endpspicture\]
where $\Delta$ denotes the diagonal.
\end{myexample}

\begin{myremark}
 Here is a useful way of thinking about this example that points us in the direction we need.  The idea is that our original pullback $BA$
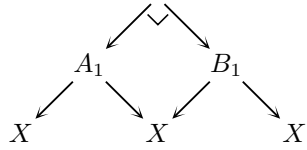
\[
\psset{unit=0.09cm,labelsep=3pt,nodesep=2pt}
\pspicture(-20,0)(40,20)



\rput(10,20){\rnode{a1}{$.$}}  
\rput(0,10){\rnode{a3}{$A_1$}}  
\rput(20,10){\rnode{a2}{$B_1$}}  
\rput(10,0){\rnode{a4}{$X$}}  
\rput(-10,0){\rnode{a5}{$X$}}  
\rput(30,0){\rnode{a6}{$X$}} 

\ncline{->}{a1}{a3} \nbput{{\scriptsize $$}} 
\ncline{->}{a1}{a2} \naput{{\scriptsize $$}} 
\ncline{->}{a3}{a4} \nbput{{\scriptsize $$}} 
\ncline{->}{a2}{a4} \naput{{\scriptsize $$}} 

\ncline{->}{a3}{a5} \nbput{{\scriptsize $$}} 
\ncline{->}{a2}{a6} \nbput{{\scriptsize $$}} 

\psline[linewidth=0.6pt]{-}(8.5,17.5)(10,16)(11.5,17.5)

\endpspicture
\]
ignored the fact that \Fop\ is in both $A$ and $B$.  So in fact we want a coequaliser 
\[
\psset{labelsep=2pt}
\pspicture(70,7)
\rput(0,3){\rnode{A}{$B \circ \F\op \circ A$}}
\rput(27,3){\rnode{B}{$B \circ A$}}
\rput(50,3){\rnode{C}{$B \otimes_{\Fop} A$}}
\psset{nodesep=3pt,arrows=->}
\ncline[offset=3pt]{A}{B}\naput{$$}
\ncline[offset=-3pt]{A}{B}\nbput{$$}
\ncline{B}{C}\naput{$$}
\endpspicture
\]
where the parallel maps are derived from 
\[\Fop \tmap{\alpha} A \makebox[0pt][l]{, and}\]
\[\Fop \tmap{\beta} B\]
respectively.  To form this coequaliser we put an equivalence relation on the morphisms of $BA$; this is encapsulated in the following definition, which is a generalisation of the definition of a factorisation system over a groupoid given in \cite{rw1}.

\end{myremark}

%
%
%
%
%
%

\begin{mydefinition}\label{factoverX}
 Let $C$ be a category, $J$ a subcategory with the same objects as $C$ (lluf).  A \demph{factorisation system over $J$} on $C$ consists of
\begin{itemize}
 \item a lluf subcategory $L$ of $C$ containing $J$, and
\item a lluf subcategory $R$ of $C$ containing $J$
\end{itemize}
such that every morphism in $C$ can be expressed as
\[\tmap{\in L} \tmap{\in R}\]
uniquely up to zigzags in $J$ as shown in the following diagram, where the morphisms on the left hand half of the diagram are all in $L$, those on the right are all in $R$, and the vertical dotted morphisms are in $J$.  The triangles on the left commute in $L$ and those on the right commute in $R$.
\[
\psset{unit=0.12cm,labelsep=2pt, linewidth=0.8pt}
\pspicture(0,0)(40,40)



\rput(-10,20){\rnode{a3}{$.$}}  
\rput(30,20){\rnode{a2}{$.$}}  

\rput(10,40){\rnode{a1}{$.$}}  
\rput(10,33.5){\rnode{b1}{$.$}}
\rput(10,27){\rnode{b2}{$.$}}
\rput(10,20.5){\rnode{b3}{$.$}}
\rput(10,14){\rnode{b4}{$.$}}
\rput(10,6.5){\rnode{b6}{$.$}}
\rput(10,0){\rnode{a4}{$.$}}  

\rput(10,11.5){$\vdots$}

\psset{nodesepA=5pt, nodesepB=1pt}

\ncline{->}{a3}{a1} \naput{{\scriptsize $\in L$}} 
\ncline{->}{a3}{b1} \naput{{\scriptsize $$}}
\ncline{->}{a3}{b2} \naput{{\scriptsize $$}}
\ncline{->}{a3}{b3} \naput{{\scriptsize $$}}
\ncline{->}{a3}{b4} \naput{{\scriptsize $$}}
\ncline{->}{a3}{b6} \naput{{\scriptsize $$}}
\ncline{->}{a3}{a4} \nbput{{\scriptsize $$}} 

\psset{nodesepA=1pt, nodesepB=5pt}

\ncline{->}{a1}{a2} \naput{{\scriptsize $\in R$}} 
\ncline{->}{b1}{a2} \naput{{\scriptsize $$}}
\ncline{->}{b2}{a2} \naput{{\scriptsize $$}}
\ncline{->}{b3}{a2} \naput{{\scriptsize $$}}
\ncline{->}{b4}{a2} \naput{{\scriptsize $$}}
\ncline{->}{b6}{a2} \naput{{\scriptsize $$}}
\ncline{->}{a4}{a2} \nbput{{\scriptsize $$}} 

\psset{linewidth=0.9pt,nodesep=2pt,linestyle=dashed, dash=2pt 1.5pt}

\ncline[]{->}{a1}{b1} \nbput{{\scriptsize $$}} 
\ncline{->}{b2}{b1} \nbput{{\scriptsize $$}} 
\ncline{->}{b2}{b3} \nbput{{\scriptsize $$}} 
\ncline{->}{b4}{b3} \nbput{{\scriptsize $$}} 
\ncline{->}{b6}{a4} \nbput{{\scriptsize $$}}

\endpspicture\]

\end{mydefinition}

\begin{myexamples} \ \\[-1em]
\begin{enumerate}

\item If $J$ is a groupoid, we get a factorisation system over $J$ as in 
\cite{rw1}. (The authors stop just short of making this definition 
although they have all the machinery in place to make it---they have other uses 
in mind and make the following construction instead.)

 \item If $J$ is the groupoid of all isomorphisms in $C$, we get an orthogonal factorisation system in the usual sense.

\item If $J$ is all identities we get a strict factorisation system.

\item Weak factorisation systems are not in general an example, for in a weak factorisation system factorisations are unique up to diagonal fillers, or ``solutions'' of certain lifting problems, but these diagonal fillers are not necessarily in $L$ or $R$; to be a factorisation system over $J$ these fillers would need to be in $J$ and hence in both $L$ and $R$.

\end{enumerate}

\end{myexamples}

\begin{mydefinition}\label{def2} (\scp{``fs''})
 Let $A$, $B$ and $C$ be \lts. Then we say $C$ is a \demph{composite} of $A$ and $B$ if $(A,B)$ forms a \fs\ over $\F\op$ on $C$. In this case we say we have a \demph{\dl} of $A$ over $B$.
\end{mydefinition}

\begin{proposition}\label{compositefs}
Given any category $C$ with a \fs\ over $\F\op$ given by $(A,B)$, if $A$ and $B$ are \lts\ then $C$ is also a \lt.
\end{proposition}

\noi As before (for the definition in $\Prof\op$), we need to check the necessary facts about finite products.  Again we defer this proof until later (Corollary~\ref{coh}).

\begin{remark}
 Note that the natural way of stating this definition involved starting with a category $C$ and ``decomposing it'' via a \fs\ over \Fop, rather than starting with \lts\ $A$ and $B$ and ``combining them'' as in other definitions.  This different viewpoint could shed light on the question of when an algebraic theory can can be expressed as a composite of simpler ones, as opposed to when it is ``irreducible''.  
\end{remark}

\noi In any case the formulation as a coequaliser gives us a good abstract formalism.  Effectively we have taken the monoidal category $\Span(\Fop, \Fop)$, put a new tensor product $\otimes_{\Fop}$ on it, and taken \dls\ with respect to this.  This is more elegantly described using bimodules.

\begin{proposition}
A \lt\
\[\Fop \tmap{\alpha} A\]
is an $(\Fop, \Fop)$-bimodule in \Span.  
\end{proposition}

\noi Then $\otimes_{\Fop}$ described above is just bimodule composition.  Thus 
the above definition of \dl\ amounts to regarding $A$ and $B$ as monads in 
$\Mod(\Span)$ and taking \dls\ between them.  But we know $\Mod(\Span) \catequiv 
\Prof\op$, so we have proved the following theorem.

\begin{theorem}\label{equiv12}
 \Dls\ as in Definition~\ref{def2} \scp{``fs''} are equivalent to \dls\ as in Definition~\ref{def1} \scp{``prof''}.
\end{theorem}

\noi We will state this more precisely later in terms of comparison functors, 
but the idea is that Definition~\ref{def2} \scp{``fs''} can be taken as an 
explicit characterisation of Definition~\ref{def1} \scp{``prof''}.  

\begin{myremarks} \ \\[-1.5em]
\begin{enumerate} 
\item This definition generalises the definition of ``\dl\ with respect to $J$'' given in \cite{rw1} although there it is expressed quite differently.  $J$ is required to be a groupoid in order to yield an equivalence relation on the morphisms of $BA$.  Effectively, this is to get unique factorisations up to plain morphisms in $J$ rather than zigzags (see \cite[Section 5.4]{rw1}).  In fact the authors do not actually mention factorisation systems over general groupoids---their aim is to give a bicategory in which orthogonal \fss\ are \dls, so once they have this general notion of \dl\ in place, they set $J$ to be the groupoid of all isomorphisms for the purposes of the \fs. 

\item Lack discusses a version of this in \cite[Sections 4.2, 4.3]{lac4}.  He is mostly concerned with PROPs, so only mentions this in passing, and again only in the case where $J$ is a groupoid.  However, his subsequent sections study \dls\ in $\ProfMon$, which is also the subject of our next section.

\end{enumerate}
\end{myremarks}

%
%
%
%


\section{Monads in monoidal profunctors}\label{profmon}

In this section we give an approach that deals a little more explicitly with 
the finite products, by taking profunctors in monoidal categories.  These are 
defined using the definition of profunctors in \cE\ (Definition~\ref{ben}) and 
taking $\cE = \Mon$, the 
category of monoids and monoid homomorphisms. Note that a 0-cell in 
$\Prof(\Mon)$ is 
an internal category in $\Mon$, that is, a strict monoidal category.

%
%
%
%
%
%
%
%

\begin{proposition}
 A monad in $\Prof(\Mon)\op$ on a monoidal category $X$ consists of a strict monoidal category $A$ and an identity-on-objects strict monoidal functor $X \tra A$.
\end{proposition}

\begin{prf}
 Follows from Theorem~\ref{theoremA}.  Put $K = \Span(\Mon)$, and $x = \fn{ob}X$, so a monoid in $K(x,x)$ in this case is a strict monoidal category with the same objects as $X$.  A morphism of such monoids is a strict monoidal identity-on-objects functor.  \end{prf}

\noi The following result is analogous to Corollary~\ref{ltmonadinprof}.

\begin{corollary}
 Every \lt\ is a monad in $\Prof(\Mon)\op$ on the 0-cell $\F\op$ regarded as a 
monoidal category with respect to product.  Conversely a monad in $\Prof(\Mon)\op$ on 
the 0-cell $\F\op$ is a strict monoidal category $A$ equipped with an 
identity-on-objects, strict monoidal functor $\F\op \tmap{\alpha_A} A$; it is a 
Lawvere theory precisely if the monoidal structure on $A$ is given by finite 
products.
 \end{corollary}

\noi Comparing this situation with that of monads in plain $\Prof\op$ we see that 
the monoidal framework is slightly ``better'': monads in 
$\Prof(\Mon)\op$ are slightly closer to being Lawvere theories in the sense that we 
only need to check a condition on $A$ and the condition on $\alpha_A$ is then 
automatic. In the next section we will give an even ``better'' framework in which 
all the conditions are automatic.

\begin{mydefinition}\label{defprofmon} (\scp{``profmon''})
 Given Lawvere theories $A$ and $B$, a \demph{distributive law} of $A$ over $B$ is a distributive law of $A$ over $B$ expressed as monads in $\ProfMon\op$. The iterated version is defined likewise, as in Theorem~\ref{idl}.  \end{mydefinition}

\begin{proposition}\label{compositeprofmon}
 The resulting composite monad $B \otimes_{\Fop} A$ is also a Lawvere theory.
\end{proposition}

\begin{myremarks} \ \\[-1.5em]

\begin{enumerate}

\item An immediate question is whether or not this gives the same thing as \dls\ in 
plain \Prof.  The perhaps surprising answer is that they are indeed the same, as 
when the monoidal structure is product, natural transformations are automatically 
monoidal.  We will discuss this in Section~\ref{comparison}.

\item  This approach is closely related to Lack's approach to 
distributive laws for PROPs in \cite{lac4}.  For PROPs, instead of 
$\F$ we use $\bP$, a skeleton of the the category of finite sets and 
\emph{bijections}.  The rest of the formalism is the same.

\end{enumerate}

\end{myremarks}

\noi As before, we defer the proof that the composite is a \lt\ until 
Section~\ref{comparison}, but it is instructive to compare the question to the 
analogous question in $\Prof\op$.  There, the issue was both whether the composite 
had finite products and whether the identity-on-objects functor preserved them.  This 
time, we know the identity-on-objects functor must preserve the monoidal structure 
of the composite, so we only need to check that this monoidal structure is given by 
finite products.

There are (at least) two ways to prove this.  A direct hands-on 
method might be possible, but a more abstract approach uses a free 
finite-product category 2-monad.  This is the subject of the next 
section.


\section{Monads in a Kleisli bicategory of profunctors}\label{kleisli}

%
%
%

In this section we follow \cite{hyl1} and use a bicategory in which monads are 
precisely Lawvere theories.  (This statement allows for types---for untyped 
Lawvere theories we will restrict to the 0-cell $1$.)

The idea in in \cite{hyl1} is to consider notions of algebraic theory determined 
by 2-monads $S$ on the 2-category \Cat\ of small categories.  If $S$ extends to 
a pseudo-monad $S_P$ on \Prof\ in a suitable way, then many-sorted 
$S$-algebraic theories arise as monads in $\Kl(S_P)$. One example is when $S$ 
is the 2-monad for strictly associative products, in which case the 
$S$-algebraic theories in this sense are (many-sorted) Lawvere theories.

Suitable extensions of $S$ to \Prof\ are given by a generalisation of 
distributive laws for monads.  The idea is that the presheaf functor sending a 
small category \bC\ to $[\bC\op, \Set]$ is almost a pseudomonad other than size 
issues, as it is in fact a pseudofunctor $\Cat \tra \CAT$, from small 
categories to locally small categories.  The bicategory \Prof\ is essentially 
the Kleisli bicategory for this not-quite monad.  

Hyland makes this precise by defining a notion of Kleisli structure on an 
inclusion of bicategories.  The idea builds from the Kleisli formulation of a 
monad given in \cite{man1}. This has the advantage of being 
applicable 
even when structure is only defined on a subcollection of objects, giving rise 
to the relative monads of \cite{acu1}.  Kleisli structures are a 
2-dimensional version of relative monads.  

The presheaf construction is a key 
example. For a small category $A$ write $PA = [A\op, \Set]$.  The following 
results are all from \cite{hyl1}.

\begin{proposition} {\bfseries (Hyland \cite{hyl1})} The presheaf construction 
$P$ gives a Kleisli structure on the inclusion $\Cat \tra \CAT$ and the 
resulting Kleisli bicategory $\Kl(P) \iso \Prof$. \end{proposition} 

\noi We will not need any details about Kleisli structures; we just need the 
following results.

\begin{proposition} {\bfseries (Hyland \cite{hyl1})}
 Let \cF\ be the monad for strictly associative products on \Cat.  This extends 
to a pseudomonad $\cF_P$ on \Prof.
\end{proposition}

\noi By abuse of notation we will also write the extended pseudomonad as \cF; 
this should not cause ambiguity as we will never need to use the original monad 
on $\Cat$. 

\begin{myremarks}\label{useful}

It is useful to take a moment to make some of the structure of \cF\ explicit; we will need this in the proof of Proposition~\ref{keyprop}.

\begin{enumerate}

\item First we make explicit the structure of $\cF A$ where $A$ is any category.  Objects in $\cF A$ are finite strings of objects in $A$.  Since these are to be products, a morphism
\[(a_1, \ldots, a_n) \tra (b_1, \ldots b_m)\]
is given by 
\begin{itemize}
\item for each index on the right a choice of projection from the left; that is a function $\alpha\: [m] \tra [n]$, and
\item for each $i \in [m]$ a morphism $a_{\alpha(i)} \tra b_i$ in $A$.
\end{itemize}
In the proof of Proposition~\ref{keyprop} we will need the morphisms of $\cF ^2 1$.  An object in this category is a string of natural numbers.  We see that in this case a morphism
\[(a_1, \ldots, a_n) \tra (b_1, \ldots b_m)\]
is given by 
\begin{itemize}
\item a function $\alpha\: [m] \tra [n]$, and
\item for each $i \in [m]$ a function $[b_i] \tra [a_{\alpha(i)}]$.
\end{itemize}

 \item Next we give the action of \cF\ on morphisms.  Given a profunctor
\[F \: A \pra B \mbox{ \ \ i.e. \ \ } B\op \times A \tra \set\]
we need a profunctor
\[\cF F \: \P A \pra \cF B \mbox{ \ \ i.e. \ \ } \cF B\op \times \cF A \tra \set.\]
The profunctor $\cF F$ is defined by 
\[\begin{array}{rcl}
\cF F(b_1, \ldots, b_n ; a_1, \ldots, a_n) &=& \dis\coprod_{\alpha \in \set(m,n)} \nd \dis\prod_{j \in [m]} F(b_{\alpha_j}, a_j) \\[18pt]
&=& \dis\prod_{j \in [m]} \nd \dis\coprod_{i \in [n]} F(b_i, a_j) \\[18pt]
&=& \dis\prod_{j \in [m]} \nd \left( \coprod_{i \in [n]} F(b_i, a_j) \right)^m 
\end{array}\]

\item Next we give the monad structure.  For multiplication we have
\[\mu \: \cF^2 1 \pra \cF 1 \mbox{ \ \ i.e. \ \ } \cF 1\op \times \cF^2 1 \tra \set\]
given by
\[\mu(n; k_1, \ldots, k_m) = \cF 1 (n, k_1 + \cdots + k_m).\]
For the unit we have
\[\eta \:  1 \pra \cF 1 \mbox{ \ \ i.e. \ \ } \cF 1\op \tra \set\]
given by
\[\eta(k) = \cF 1(k,1) = \set(1,[k]).\]

\end{enumerate}

\end{myremarks}

\begin{mydefinition}
From henceforth we shall write \ProfP\ for $\Kl(\cF_P)$, the Kleisli 
bicategory of \cF\ extended to \Prof.
\end{mydefinition}

\noi Monads in \Profp\ are then many-sorted Lawvere theories; we only need the 
following special case.

\begin{theorem}{\bfseries (Hyland)} \label{hylandprofp}
 Monads on 1 in \Profp\ are 
precisely un-typed \lts.
\end{theorem}


\begin{prf} {\bfseries (Sketch.)} A 1-cell $A \tra B$ in $\Kl(\cF_P)$ 
is a profunctor $A \pra \cF B$, i.e. a functor $\cF B\op \times A \tra \Set$.  
So a monad on $1$ has an underlying functor $\cF 1\op \times 1 \tra \set$, i.e. 
a functor $\FinSet \tra \Set$ or equivalently a finitary functor $\Set \tra 
\Set$; the monad structure then makes this into a finitary monad on \Set. 
\end{prf}

%
%

\noi In fact we have a more precise result involving an equivalence of categories (Theorem~\ref{thetaequiv}).  Before we prove that, the following proposition provides a functor that will evaluate a monad in $\Profp\op$ at the corresponding \lt\ expressed in $\Prof\op$.  Recall that the forgetful functor from the Kleisli category of any monad to its underlying category is given on morphisms by applying the monad and postcomposing with $\mu$.  The following proposition evaluates this for $\cF$. 

\begin{proposition}\label{keyprop}
 
For any profunctor $1 \pmap{F} \cF 1$, the composite
\[ \cF 1 \pmap{\cF  F} \cF ^2 1 \pmap{\mu} \cF 1\]
 is given by 
\[\begin{array}{ccc}
   \cF 1\op \times \cF 1 &\tra& \set\\
(j,n) &\tmapsto& \set(n,Fj)
  \end{array}.\]
\end{proposition}

\begin{prf}
 By definition this composite is

\[\begin{array}{rcl}
   (j,l) & \tmapsto &  \dis\int\limits^{(k_1, \ldots, k_m) \in \cF^21} \hh{-2em} \mu(j; k_1 + \cdots +k_m) \times \cF F(k_1, \ldots, k_m; l)\\[12pt]
&=& \dis\int\limits^{(k_1, \ldots, k_m) \in \cF^2 1} \hh{-2em} \cF 1(j, k_1 + \cdots + k_m) \times \left( \coprod_{i \in [m]} F(k_i) \right)^l
  \end{array}\]

We aim to show that in computing this coend we only need to consider $m=1$.  We use the fact that in general in a coend cocone for $Q: \bI\op \times \bI \tra \set$ 

\[
\psset{unit=0.1cm,labelsep=1pt,nodesep=2pt,npos=0.4}
\pspicture(0,0)(40,20)



\rput(20,20){\rnode{a1}{$Q(U,U)$}}  
\rput(0,10){\rnode{a3}{$Q(V,U)$}}  
\rput(40,10){\rnode{a2}{$.$}}  
\rput(20,0){\rnode{a4}{$Q(V,V)$}}  

\ncline{->}{a3}{a1} \naput{{\scriptsize $Q(f,1)$}} 
\ncline{->}{a3}{a4} \nbput{{\scriptsize $Q(1,f)$}} 
\ncline{->}{a1}{a2} \naput{{\scriptsize $$}} 
\ncline{->}{a4}{a2} \nbput{{\scriptsize $$}} 

\endpspicture
\]
for $f \: U \tra V$ in $\bI$, if $Q(1,f)$ is surjective we can ignore $Q(V,V)$ as no further information is contributed by it.

Choose $U,V \in \cF^2 1$ as follows
\[\begin{array}{rcl}
U &=& (k_1 + \cdots + k_m) = (k), \mbox{ \ say} \\
V &=& (k_1, \ldots, k_m).
     \end{array}\]
Note that $U$ is a 1-ary string.  We then define $f \: U \tra V \in \cF^2 1$ as follows.  Recall that a morphism 
\[(a_1, \ldots, a_n) \tra (b_1, \ldots, b_m)\]
in $\cF ^21$ consists of  
\begin{itemize}
 \item a map $\alpha\: m \tra n$ in \set, and 
\item for all $i \in [m]$, a map $\beta_i \: b_i \tra a_{\alpha(i)}$ in \set. 
\end{itemize}
Here we have $n=1$ so $\alpha$ is trivial, thus to define $f$ we just need to give, for all $i \in [m]$ a map
\[\beta_i \: k_i \tra k_1 + \cdots + k_m \in \set\]
and we set these to be the canonical coproduct insertions.

Now note that
\[ \begin{array}{rcl}
    Q(V,U) &=& \cF 1(j, k_1 + \cdots + km) \times \left( \dis\coprod_{i \in [m]} F(k_i) \right)^l\\
&\iso & Q(V,V) 
   \end{array}\]
and moreover the isomorphism is given by $Q(1,f)$.  So we can disregard all vertices in the coend for which $m \neq 1$.

Thus the coend becomes
\[\int\limits^{k \in \cF 1} \hh{-0.5em} \cF 1(j,k) \times \big(F(k)\big)^l \iso \set(n,Fj)\]
as required.\end{prf}

\begin{remark}
Note that this profunctor will be called $\bar{F}$ in Section~\ref{comparison} and it will give us the comparison between the profunctor approach and the monad approach; note that if $F$ is a finitary monad, $\bar{F}$ is its associated \lt.
\end{remark}

\noi Write $[\Set,\Set]_f$ for the monoidal category of 
finitary endofunctors on \Set\ and natural transformations, with the monoidal 
structure given by composition.

\begin{theorem}\label{equivcatprof}\label{thetaequiv}
 There is a monoidal equivalence of categories
\[ [\set,\set]_f \tmap{\catequiv} \profp\op(1,1).\]
\end{theorem}

\begin{prf}
Recall that a finitary functor $F\: \set \tra \set$ is entirely determined by its restriction to \finset, by the formula
\[FX = \int\limits^{[n] \in \finset} \hh{-1.6em} F[n] \times X^n.\]
We define a functor
\[ [\set,\set]_f \tmap{\theta} \profp\op(1,1)\]
as follows. Given a finitary functor $F: \set \tra \set$ we restrict it as
\[\cF  1 \op \catequiv \finset \tmap{F} \set \]
which can be regarded as a profunctor $1 \pmap{\theta F} \cF  1$ as required.  (Note 
that technically we must pick a functor $\cF  1 \op \tra \finset$ giving the 
equivalence.) On morphisms we also take the restriction of natural transformations 
to \finset. 

The interesting part is the monoidal structure, which is given by composition.  Consider finitary functors 
\[\set \tmap{F} \set \tmap{G} \set.\]
Then the composite $\theta{G} \circ \theta{F}$ in $\profp$ is given by the profunctor composite
\[1 \pmap{\theta{F}} \cF 1 \pmap{\cF (\theta{G})} \cF ^2 1 \pmap{\mu} \cF 1\]
which is some functor
\[\cF  1\op \tra \set.\]
Now, using the formula for $\mu$ and the action of \cF \ on morphisms as given in Remarks~\ref{useful} we see that the composite is given by
\[\begin{array}{rcll}
m &\tmapsto& \dis\int\limits^{j \in \cF  1} \set(j, \theta{G}(m)) \times \theta{F}(j) &\mbox{by Proposition~\ref{keyprop}}\\[12pt]
   &=& \dis\int\limits^{j \in \cF  1} \set(j, Gm) \times F(j) \\[12pt]
&=& FG(m) & \mbox{by standard density}\\[6pt]
 &=& \theta(FG)(m)
  \end{array}\]

\noi Full and faithfulness is clear; essential surjectivity of $\theta$ follows from 
the fact that a finitary functor $F$ is determined uniquely up to isomorphism by its 
restriction to \finset. \end{prf}

\begin{definition}\label{def3} \scp{``kleisli''}
 Given \lts\ $A$ and $B$, a \demph{\dl\ of $A$ over $B$} is a \dl\ of $A$ over $B$ 
expressed as monads on 1 in $\profp\op$ via Theorem~\ref{hylandprofp}.  The 
composite monad $BA$ is automatically a \lt, and is called the \demph{composite \lt}. 
The iterated version is defined likewise, as in Theorem~\ref{idl}.
\end{definition}

\noi Note that this is the only case in which it is immediate that the composite 
monad is a \lt; however the result for the other definitions will follow.  First, we 
can immediately deduce from the preceding results that these \dls\ correspond 
precisely to \dls\ between finitary monads in \set.

\begin{corollary}\label{equiv3monad}
 Let $S$ and $T$ be finitary monads on \set\ with associated \lts\
\[\begin{array}{rcl}
   \theta(S) &=& \L_S\\
\theta(T) &=& \L_T
  \end{array}\]
expressed as monads on 1 in $\profp\op$.  Let
\[\lambda \: ST \trta TS\]
be a \dl\ of $S$ over $T$.  Then
\[\theta(\lambda) \: \theta(ST) \trta \theta(TS)\]
gives a \dl\ of $\L_S$ over $\L_T$ in $\profp\op$ as
\[\begin{array}{rcl}
   \theta(ST) &\cong& \L_S \L_T \\
\theta(TS) &\cong& \L_T \L_S.
  \end{array}\]
Furthermore since $\L_{TS} = \theta(TS) \cong \L_T \L_S$ we see that the composite \lt\ is the \lt\ associated to the composite monad.  Conversely since $\theta$ is an equivalence, every \dl\ of \lts\ arises in this way. 
\end{corollary}

\begin{myremark}
 In fact since \dls\ in a 2-category $K$ are the 0-cells of $\Mnd(\fn{Mnd} K)$ we 
could express this as a biequivalence between the ``bicategories of \dls'', and then 
iterate the construction to get a notion of iterated \dl\ for \lt, as in 
Definition~\ref{idl}.  
\end{myremark}


\section{Comparison}\label{comparison}

We now have four definitions of \dl\ for \lt\ in place:
\begin{enumerate}
 \item \dprof: \Dls\ in $\Prof\op$.
\item \dfs: \Fss\ over $\F\op$.
\item \dprofmon: \Dls\ in $\ProfMon\op$.
\item \dkl: \Dls\ in $\profp\op$.
\end{enumerate}
So far we have shown that 
\begin{itemize}
 \item \dprof\ and \dfs\ are equivalent (Theorem~\ref{equiv12}).
\item \dkl\ is equivalent to the monad approach (Corollary~\ref{equiv3monad}).
\end{itemize}
In this section we will complete the programme of equivalences by showing that \dprof\ is equivalent to both \dprofmon\ and the monad approach.  The following diagram shows comparison functors we will construct; so far we have exhibited $\theta$: 

\[
\psset{unit=0.12cm,labelsep=2pt,nodesep=3pt, linewidth=0.8pt}
\pspicture(0,-20)(52,25)


\rput(-15,1){\parbox{6em}{\begin{center}``monad \ \hh{2pt}
 \\approach''\end{center}}}

\rput(3,1){\rnode{a3}{$[\Set,\Set]_f$}}

\rput(40,25){\rnode{a2}{$\Profp\op(1,1)$}}   
\rput(40,8){\rnode{b2}{$\Prof\op(\cF 1, \cF 1)$}}

\rput[l](60,25.5){\scp{``kl''}}
\rput[l](60,1.5){\scp{``prof''}}
\rput[l](60,-17.5){\scp{``profmon''}}

\rput(37,5){\scalebox{0.8}{\rotatebox{270}{$\catequiv$}}}
 
\rput(40,1){\rnode{a4}{$\Prof\op(\Fop,\Fop)$}}  
\rput(40,-18){\rnode{b1}{$\ProfMon\op(\Fop,\Fop)$}}  

\ncline{->}{a3}{a2} \naput{{\scriptsize $\theta$ \textsf{equivalence}}}  
\ncline{->}{a3}{a4}\naput{{\scriptsize $\phi$ \textsf{f+f}}}
\ncline{->}{a2}{b2} \naput{{\scriptsize $U^\theta$ \textsf{{forgetful}}}}
\ncline{->}{b1}{a4}\nbput{{\scriptsize $U^\psi$ \textsf{{forgetful}}}} 

\ncline{->}{a3}{b1} \nbput{{\scriptsize $\psi$ \textsf{f+f}}}

\endpspicture
\]

First we make explicit the functor $\phi$ as follows.
\[\begin{array}{ccc}
   \phi \: [\Set, \Set]_f & \tra & \Prof\op(\F\op, \F\op) \\
F & \tmapsto & \bar{F} \: \F \times \F\op \tra \Set \\
&& \bar{F}(n,m) = \Set(m, Fn) \\
\alpha:F \Tra G & \tmapsto & 
\psset{unit=0.1cm,labelsep=2pt,nodesep=3pt}
\pspicture(-5,10)(20,20)
\rput(0,10){\rnode{a1}{$\F \times \F\op$}}  
\rput(20,10){\rnode{a2}{$\Set$}}  
\pcline[doubleline=true,arrowinset=0.6,arrowlength=0.7,arrowsize=0.8pt 2.6]{->}(10,14)(10,6)\naput{{\scriptsize $\bar{\alpha}$}}

\ncarc[arcangle=45]{->}{a1}{a2}\naput{{\scriptsize $\bar{F}$}}
\ncarc[arcangle=-45]{->}{a1}{a2}\nbput{{\scriptsize $\bar{G}$}}

\endpspicture\\[24pt]
\alpha_n \: Fn \tra Gn && \bar{\alpha}_{n,m} \: \bar{F}(n,m) \tra \bar{G}(n,m) \\
&& \Set(m, Fn) \ltmap{\alpha_n \circ \uscore} \Set(m, Gn)
  \end{array}\]

\noi Later we will show that this is a monoidal functor, but now we concentrate on 
other properties.

\begin{proposition} \label{phifull}
 The functor $\phi$ is clearly faithful (by Yoneda).  It is also full.
\end{proposition}


\begin{prf}
 Suppose we have a natural transformation
\[\psset{unit=0.1cm,labelsep=2pt,nodesep=3pt}
\pspicture(20,20)
\rput(0,10){\rnode{a1}{$\F \times \F\op$}}  
\rput(20,10){\rnode{a2}{$\Set.$}}  
\pcline[doubleline=true,arrowinset=0.6,arrowlength=0.7,arrowsize=0.8pt 2.6]{->}(10,14)(10,6)\naput{{\scriptsize ${\beta}$}}

\ncarc[arcangle=45]{->}{a1}{a2}\naput{{\scriptsize $\bar{F}$}}
\ncarc[arcangle=-45]{->}{a1}{a2}\nbput{{\scriptsize $\bar{G}$}}

\endpspicture\]

\noi We aim to show that $\beta$ is in fact of the form $\bar{\alpha}$ for some $\alpha$ as above.  Now, given any $n \in \F$ we define putative components $\alpha_n$ to be the components $\beta_{n,1}$ as shown.
\[\begin{array}{ccccc}
   \beta_{n,1} & : & \bar{F}(n,1) & \tra & \bar{G}(n,1) \\
\mbox{ie} && \Set(1, Fn) & \tra & \Set(1,Gn) \\
&& \rotatebox{90}{=} && \rotatebox{90}{=} \\
&& Fn && Gn \\
  \end{array}\]

\noi We claim

\begin{enumerate}
 \item These $\alpha_n$ are components of a natural transformation $\alpha \: F \Tra G$, and
\item $\beta = \bar{\alpha}$.
\end{enumerate}

\noi For the first part we need to check that for all $f \: n \tra k \in \F$ the following naturality square commutes
\[
\psset{unit=0.1cm,labelsep=3pt,nodesep=3pt}
\pspicture(40,18)



\rput(0,18){\rnode{a1}{$Fn$}}  
\rput(20,18){\rnode{a2}{$Gn$}}  
\rput(0,0){\rnode{a3}{$Fk$}}  
\rput(20,0){\rnode{a4}{$Gk$}}  

\ncline{->}{a1}{a2} \naput{{\scriptsize $\alpha_n$}} 
\ncline{->}{a3}{a4} \nbput{{\scriptsize $\alpha_k$}} 
\ncline{->}{a1}{a3} \nbput{{\scriptsize $Ff$}} 
\ncline{->}{a2}{a4} \naput{{\scriptsize $Gf$}} 

\endpspicture
\]

\noi Now by naturality of $\beta$ we have
\[
\psset{unit=0.1cm,labelsep=3pt,nodesep=3pt}
\pspicture(40,22)



\rput(0,20){\rnode{a1}{$\Set(Fn, Fn)$}}  
\rput(35,20){\rnode{a2}{$\Set(Fn, Gn)$}}  
\rput(0,0){\rnode{a3}{$\Set(Fn, Fk)$}}  
\rput(35,0){\rnode{a4}{$\Set(Fn,Gk)$}}  

\ncline{->}{a1}{a2} \naput{{\scriptsize $\alpha_n \circ \uscore$}} 
\ncline{->}{a3}{a4} \nbput{{\scriptsize $\alpha_k \circ \uscore$}} 
\ncline{->}{a1}{a3} \nbput{{\scriptsize $Ff \circ \uscore$}} 
\ncline{->}{a2}{a4} \naput{{\scriptsize $Gf \circ \uscore$}} 

\endpspicture
\]

\noi so starting with the identity in the top left we have
\[
\psset{unit=0.1cm,labelsep=3pt,nodesep=3pt}
\pspicture(40,22)

%


\rput(0,21){\rnode{a1}{$1_{Fn}$}}  
\rput(35,21){\rnode{a2}{$\alpha_n$}}  
\rput(0,0){\rnode{a3}{$Ff$}}  

\rput(35,6){\rnode{a5}{$Gf \circ \alpha_n$}}  

\rput(35,3){\rotatebox{90}{$=$}}

\rput(35,0){\rnode{a4}{$\alpha_k \circ Ff.$}}  

\ncline{|->}{a1}{a2} \naput{{\scriptsize $$}} 
\ncline{|->}{a3}{a4} \nbput{{\scriptsize $$}} 
\ncline{|->}{a1}{a3} \nbput{{\scriptsize $$}} 
\ncline{|->}{a2}{a5} \naput{{\scriptsize $$}} 

\endpspicture
\]

\noi Now we need to show that
\[\beta_{n,m} \: \bar{F}(n,m) \tra \bar{G}(n,m) \]
is $\alpha_n \circ \uscore$, that is, $\beta_{n,1} \circ \uscore$.  Now given $f\: m \tra Fn$ and $i \: 1 \tra m$ we have

\[
\psset{unit=0.1cm,labelsep=3pt,nodesep=3pt}
\pspicture(40,25)


\rput(0,20){\rnode{a1}{$\Set(m, Fn)$}}  
\rput(35,20){\rnode{a2}{$\Set(m, Gn)$}}  
\rput(0,0){\rnode{a3}{$\Set(1, Fn)$}}  
\rput(35,0){\rnode{a4}{$\Set(1,Gn)$}}  

\ncline{->}{a1}{a2} \naput{{\scriptsize $\beta_{m,n} $}} 
\ncline{->}{a3}{a4} \nbput{{\scriptsize $\beta_{n,1} = \alpha_n \circ \uscore$}} 
\ncline{->}{a1}{a3} \nbput{{\scriptsize $\uscore \circ i$}} 
\ncline{->}{a2}{a4} \naput{{\scriptsize $\uscore \circ i$}} 

\endpspicture
\]

\[
\psset{unit=0.1cm,labelsep=3pt,nodesep=3pt}
\pspicture(40,25)

%


\rput(0,21){\rnode{a1}{$f$}}  
\rput(35,21){\rnode{a2}{$\beta_{m,n}(f)$}}  
\rput(0,0){\rnode{a3}{$f \circ i$}}  

\rput(35,6){\rnode{a5}{$(\beta_{m,n}(f))(i)$}}  

\rput(35,3){\rotatebox{90}{$=$}}

\rput(35,0){\rnode{a4}{$(\alpha_n \circ f)(i).$}}  

\ncline{|->}{a1}{a2} \naput{{\scriptsize $$}} 
\ncline{|->}{a3}{a4} \nbput{{\scriptsize $$}} 
\ncline{|->}{a1}{a3} \nbput{{\scriptsize $$}} 
\ncline{|->}{a2}{a5} \naput{{\scriptsize $$}} 

\endpspicture
\]

\noi This is true for all $i \in m$, so $\beta_{m,n}(f)$ and $\alpha_n \circ f$ 
agree everywhere, hence $\beta = \bar{\alpha}$ as required and the functor $\phi$ is 
indeed full. \end{prf}

\begin{proposition}
 The functor $\phi: [\Set,\Set]_f \tra \Prof\op(\F\op, \F\op)$ factors 
through 
$\ProfMon\op(\F\op, \F\op)$, giving the functor we called $\psi$ above.
\end{proposition}

\begin{prf}  We use the definition of $\ProfMon\op$ as $\Mod(\Span(\Mon))$, and 
$\Prof\op$ as 
$\Mod(\Span)$. We write the underlying span of $\F\op$ as
%
%
\[
\psset{unit=0.1cm,labelsep=3pt,nodesep=1.5pt}
\pspicture(30,15)


\rput(11,11){\rnode{a1}{$\F_1$}}  
\rput(0,0){\rnode{a3}{$\F_0$}}  
\rput(22,0){\rnode{a2}{$\F_0$}}  

\ncline{->}{a1}{a3} \nbput{{\scriptsize $t$}} 
\ncline{->}{a1}{a2} \naput{{\scriptsize $s$}} 
\endpspicture\]
Now consider a finitary functor $F: \Set \tra \Set$. Now the image of $F$ under 
$\phi$ is $\bar{F}$, whose underlying span of $\bar{F}$ as an $\F\op$-bimodule is

\[
\psset{unit=0.1cm,labelsep=3pt,nodesep=1.5pt}
\pspicture(40,15)


\rput(13,13){\rnode{a1}{\scalebox{0.8}{$A= \displaystyle\coprod_{m,n}\set(m,Fn)$}}}  
\rput(0,0){\rnode{a3}{{\scriptsize $\F_0$}}}  
\rput(26,0){\rnode{a2}{{\scriptsize $\F_0$}}}  

\ncline{->}{a1}{a3} \nbput{{\scriptsize $$}} 
\ncline{->}{a1}{a2} \naput{{\scriptsize $$}} 
\endpspicture\]

\noi The claim is that this is automatically a bimodule in \Mod(\Span), 
although it is 
a priori just a bimodule in \Span.  So we need to put a monoid structure on $A$ such 
that the left and right $\F$-actions respect this.  Note that the monoid structure 
in $\F_0$ is given by addition.  So given 
\[\begin{array}{cc}
f_1 \: m_1 \tra Fn_1 \\
f_2 \: m_2 \tra Fn_2    
  \end{array}\]
we construct a function
\[f_1 \oplus f_2 \: m_1 + m_2 \tra F(n_1 + n_2).\]

Now by coproduct in \Set\ we certainly have
\[m_1 + m_2 \ltmap{f_1+f_2} Fn_1 + Fn_2 \ltmap{\mbox{{\scriptsize {\sf canonical}}}} F(n_1+n_2)\]
and we call this $f_1 \oplus f_2$.  We also need $e \: 0 \tra F0$ such that
\[f \oplus 0 = 0 \oplus f = f.\]  
This is the unique map. Then $f \oplus 0$ is the following map:
\[m = m+ 0 \ltmap{f + !} Fn + F0 \ltmap F(n+0) = Fn\]
which is the same as $f$ by a straightforward diagram chase.


Now we must check actions.  These are given by pre- and post-composition.  For the left action, given $k \tmap{f} m$ in \set\  we have 
\[\Set(m, Fn) \ltmap{\uscore \circ f} \Set(k, Fn)\]
and we need to check that 
\[(g_1 \oplus g_2) \circ (f_1 + f_2) = (g_1 \circ f_1) \oplus (g_2 \circ f_2).\]
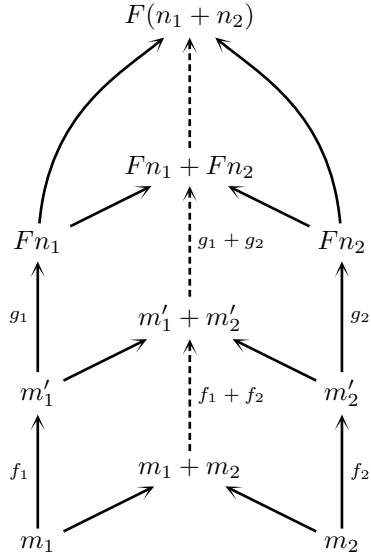
\[
\psset{unit=0.1cm,labelsep=2pt,nodesep=3pt,linewidth=1pt}
\pspicture(40,75)


\rput(0,0){\rnode{a1}{$m_1$}}  
\rput(0,20){\rnode{a2}{$m_1'$}}  
\rput(0,40){\rnode{a3}{$Fn_1$}}  

\rput(20,10){\rnode{b1}{$m_1+m_2$}}  
\rput(20,30){\rnode{b2}{$m_1'+m_2'$}}  
\rput(20,50){\rnode{b3}{$Fn_1+Fn_2$}}  
\rput(20,70){\rnode{b4}{$F(n_1+n_2)$}}  

\rput(40,0){\rnode{c1}{$m_2$}}  
\rput(40,20){\rnode{c2}{$m_2'$}}  
\rput(40,40){\rnode{c3}{$Fn_2$}}  

\ncline{->}{a1}{a2} \naput{{\scriptsize $f_1$}}
\ncline{->}{a2}{a3} \naput{{\scriptsize $g_1$}}
\nccurve[angleA=85,angleB=220,ncurvA=1,ncurvB=0.5]{->}{a3}{b4}\naput{{\scriptsize $$}}

\ncline{->}{c1}{c2} \nbput{{\scriptsize $f_2$}}
\ncline{->}{c2}{c3} \nbput{{\scriptsize $g_2$}}
\nccurve[angleA=95,angleB=320,ncurvA=1,ncurvB=0.5]{->}{c3}{b4}\naput{{\scriptsize $$}}

\ncline{->}{a1}{b1} \naput{{\scriptsize $$}}
\ncline{->}{a2}{b2} \naput{{\scriptsize $$}}
\ncline{->}{a3}{b3} \naput{{\scriptsize $$}}

\ncline{->}{c1}{b1} \naput{{\scriptsize $$}}
\ncline{->}{c2}{b2} \naput{{\scriptsize $$}}
\ncline{->}{c3}{b3} \naput{{\scriptsize $$}}

\psset{linestyle=dashed, dash=2.5pt 1.5pt}
\ncline{->}{b1}{b2} \nbput{{\scriptsize $f_1+f_2$}}
\ncline{->}{b2}{b3} \nbput{{\scriptsize $g_1+g_2$}}
\ncline{->}{b3}{b4} \nbput{{\scriptsize $$}}

\endpspicture
\]
The left and right hand sides of the equation we want then just correspond to the middle dotted composite of this diagram associated either way round, so the result follows by associativity.

For the right action, given $n \tmap{f} k$ in \set\ we have
\[\Set(m, Fn) \ltmap{Ff \circ \uscore} \Set(m, Fk)\]
and we need to check that
\[F(f_1 + f_2) \circ (g_1 \oplus g_2) = (Ff_1 \circ g_1) \oplus (Ff_2 \circ g_2).\]
The result then follows by a straightforward diagram chase involving diagrams similar 
to the previous one. This completes the result on objects. 

We must now check the result on morphisms, that is, given a natural transformation 
between finitary functors $\alpha: F \Tra G$ we must show that $\bar{\alpha}$ is a 
monoid map as
\[\dis\coprod_{m,n} \Set(m, Fn) \tra \dis\coprod_{m,n} \Set(m, Gn).\]
So we need to show that given
\[\begin{array}{c}
   f_1 \: m_1 \tra Fn_1 \\
f_2 \: m_2 \tra Fn_2
  \end{array}\]
we have
\[\alpha_{n_1 + n_2} \circ (f_1 \oplus f_2) = (\alpha_{n_1} \circ f_1) \oplus (\alpha_{n_2} \circ f_2).\]
This follows from a straightforward diagram chase and naturality of $\alpha$.  
\end{prf}

\begin{corollary}
 It follows immediately that
\[[\Set, \Set]_f \tmap{\psi} \ProfMon\op(\F\op, \F\op)\]
is full as well as faithful, and likewise the forgetful functor
\[\ProfMon\op(\F\op, \F\op) \map{U^\psi} \Prof\op(\F\op, \F\op)\]
is also full and faithful on the image of $\psi$. Thus distributive laws according to \scp{``profmon''} correspond to distributive laws according to \scp{``prof''}.
\end{corollary}

\begin{prf}
 Follows from $\phi$ being full (Proposition~\ref{phifull}).
\end{prf}

\begin{corollary}\ \label{obvious}\\[-1.5em]

\begin{enumerate}

\item  When $F$ is a monoid in $[\Set, \Set]_f$ (i.e. a finitary monad on 
\Set), $\psi{F}$ is a monad in $\ProfMon\op$ given by an identity-on-objects 
functor 
\[\alpha \: \F\op \tra A\]
where the monoidal structure of $A$ is given by products.  Conversely any such monad in $\ProfMon\op$ arises in this way.

\item  When $F$ is a monoid in $[\Set, \Set]_f$, $\phi F$ is a monad in 
$\Prof\op$ given by an identity-on-objects functor 
\[\alpha \: \F\op \tra A\]
where $A$ has finite products and $\alpha$ preserves finite products.  Conversely any such monad in $\Prof\op$ arises in this way.

\end{enumerate}

\end{corollary}

\begin{remark}
 Note that this is not much more than the standard correspondence between finitary monads and \lts.
\end{remark}

\begin{prf} Regarding $\phi F = \bar{F}$ as a category, we have $\bar{F}(n,m) = 
\Set(m,Fn)$.  We need to check that $n+k$ is the categorical product in $\bar{F}$, 
that is, there is a natural isomorphism
\[\bar{F}(p,n) \times \bar{F}(p,k) \iso \bar{F}(p, n+k)\]
that is
\[\Set(n, Fp) \times \Set(k, Fp) \iso \Set(n+k, Fp)\]
which is true by definition of coproduct in \Set. This proves both parts.
\end{prf}

That deals with the bottom half of the comparison diagram.  We now deal with the top half. Proposition~\ref{keyprop} shows that the canonical Kleisli forgetful functor $U^\theta$ makes the following triangle commute (up to isomorphism) 
\[
\psset{unit=0.1cm,labelsep=2pt,nodesep=3pt,npos=0.4}
\pspicture(0,-1)(50,17)

%

\rput(0,0){\rnode{a1}{$[\set,\set]_f$}}  
\rput(47,0){\rnode{a2}{$\Prof\op(\Fop,\Fop)$}}  
\rput(47,15){\rnode{a3}{$\Profp\op(1,1)$}}  

\ncline{->}{a1}{a2} \nbput[npos=0.55]{{\scriptsize $\phi$}} 
\ncline{->}{a1}{a3} \naput[npos=0.51]{{\scriptsize $\theta$}} 
\ncline[linestyle=dashed, dash=2.5pt 1.5pt,labelsep=3pt]{->}{a3}{a2} \naput[labelsep=2pt]{{\scriptsize $U^\theta$}} 

\endpspicture
\]
  
\begin{proposition} Distributive laws according to \scp{``prof''} correspond to those according to the monad approach. \end{proposition}

\begin{prf}
As the forgetful functor $U^\theta$ is a 1-object restriction of a 
pseudo-functor, it must be monoidal.  Thus the functor $\phi$ must be monoidal, and we have already shown that it is full and faithful.  Thus distributive laws according to \dprof\ are equivalent to those according to the monad approach.
\end{prf}

\begin{remark}
Since $\theta$ is an equivalence and $\phi$ is full and faithful, the 
forgetful functor
must also be full and faithful, giving a direct comparison between the \scp{kleisli} and \scp{prof} approaches.  
\end{remark}

This completes the suite of equivalences.

\begin{corollary}\label{coh}\label{compprof} Let $A$ and $B$ be \lts\ 
expressed according to any of Definitions~\ref{def1} ``\scp{prof}'', 
\ref{def2} ``\scp{fs}'' or \ref{defprofmon} ``\scp{profmon}'', and let 
$\sigma \: AB \tra BA$ be a \dl\ of $A$ over $B$ according to the same 
definition.  Then the composite monad $BA$ is also a \lt.
\end{corollary}

%
%

\begin{prf} For monads in $\Prof\op$ we know can write $A$ and $B$ as $\phi S$ and 
$\phi{T}$ for some finitary monads $S$ and $T$ by Corollary~\ref{obvious}.  Then by 
fullness of $\phi$ the 2-cell $\sigma$ giving the distributive law 
must be of the form $\phi(\lambda)$ for some natural transformation $\lambda: ST 
\Tra TS$; by faithfulness the axioms for $\lambda$ to be a distributive law follow 
from those for $\sigma$.  Thus 
the composite $BA$ is isomorphic to $\phi(TS)$ thus is a \lt.  The 
result for 
factorisation systems immediately follows, and that for monads in $\ProfMon\op$ 
follows in the same way \end{prf}

Although we have now completed the equivalences, we include one further
characterisation as we find it illuminating.  It is well known that there are 
two canonical identity-on-objects pseudofunctors relating \Cat\ and \Prof.
Given a functor $C \tmap{F} D$ in $\CAT$ the two functors act as follows.

\begin{enumerate}
 \item Covariant: $C \pmap{F_*} D$ in \Prof\ defined by $F_*(d,c) = D(d,Fc)$.  
This is the canonical free pseudofunctor if we regard $\Prof$ as the Kleisli 
bicategory for the presheaf Kleisli structure.

\item Contravariant: $D \pmap{F^*} C$ defined by $F^*(c,d) = F(Fc, d)$.

\end{enumerate}

\noi Thus given a monad $\set \tmap{T} \set$ we get a monad $\set \pmap{T_*} 
\set$ in \PROF\ and this could be regarded as an algebraic theory typed in \set. 
 However if we have a finitary monad we can restrict our types to the small 
category \F\ via a chosen embedding
\[\F \tmap{I} \set.\]
Then we can define a functor
\[\begin{array}{ccc}
 [\set,\set]_f & \tra & \Prof(\F, \F)\\
\set \tmap{F} \set & \tmapsto & \F \pmap{I_*} \set \pmap{F_*} \set \pmap{I^*} \F.
\end{array}\]

The following proposition shows that on monads this gives us the (opposite of) the associated \lt.

\begin{proposition}
 The above composite gives the profunctor
\[(k,n) \tmapsto \set(k,Fn).\]
\end{proposition}

\begin{prf}
This is a routine coend calculation, using the fact that $F_* I_* = (FI)_*$:
\[\begin{array}{rrll}
   (k,n) &\mapsto& \dis\int\limits^{X \in \set} I^*(k,X) \times (FI)_*(X,n) \\[12pt]
&=& \dis\int\limits^{X \in \set} \set(k,X) \times \set(X,Fn) \\[12pt]
&=& \set(k,Fn) & \mbox{by density.}
  \end{array}\]
Finally we can regard this as $\Fop \pra \Fop$ by taking it to be in $\Prof\op$ via the standard duality.
 \end{prf}

Hence we have directly constructed the functor
\[[\set,\set]_f \tra \Prof\op(\Fop, \Fop)\]
constructed previously as the composite via $\profp$, and this gives another explanation of the (slightly annoying) presence of the ``op''.  

Furthermore, that this functor is monoidal follows neatly from the finitary conditions as follows.  We need to check that, given finitary functors
\[\set \tmap{F} \set \tmap{G} \set\]              
the composite in \Prof\
\[\F \pmap{I_*} \set \pmap{F_*}\set  \pmap{G_*} \set \pmap{I^*} \F\]
is isomorphic to
\[\F \pmap{I_*} \set \pmap{F_*}\set \pmap{I^*} \F \pmap{I_*} \set \pmap{G_*} \set \pmap{I^*} \F.\]
 
In fact $G$ being finitary gives us that
\[\set  \pmap{G_*} \set \pmap{I^*} \F\]
is isomorphic to
\[\set \pmap{I^*} \F \pmap{I_*} \set \pmap{G_*} \set \pmap{I^*} \F.\]
We simply calculate the coends.  The first gives
\[\begin{array}{rrll}
   (k,X) & \mapsto & \dis\int\limits^{A \in \set} G_*(A,X) \times I^*(k,A) \\[12pt]
&=& \dis\int\limits^{A \in \set} \set(A,GX) \times \set(k,A) \\[12pt]
&=& \set(k,GX) & \mbox{by density.}
  \end{array}\]
For the second composite we have
\[\begin{array}{rrll}
   (k,X) & \mapsto & \dis\int\limits^{n \in \F} \set(k,Gn) \times \set(n,X) \\[12pt]
&=& \set(k,GX)
  \end{array}\]
as $G$ is finitary.

\section{Future work}\label{future}

This new theory of \dls\ for \lts, with its four different viewpoints, opens up various possibilities for further study.  We conclude by briefly mentioning a few.  Some work in this direction has been undertaken in \cite{bsz1}.

\begin{itemize}

\item We could seek more concrete ways of expressing distributive laws over $\Fop$ using the (quite special) properties of $\Fop$.  We could seek ``canonical forms'' for operations in the composite theory.

 \item We could study the question of when an algebraic theory can be decomposed into simpler ones, and when it is ``irreducible'',

\item We could further study iterated distibutive laws in the context of \lts.

\item We could extend the theory to any of the generalised versions of \lt.

\end{itemize}

\subsubsection*{Acknowledgements}


This work was launched by a question posed to me by Jean B\'enabou at the 89th PSSL in Louvain-la-Neuve, for which I am grateful.  Its progress was then dramatically catalysed by invitations I received to speak at the 4th Scottish Category Seminar and at ``Category Theory, Algebra and Geometry'' in Louvain-la-Neuve in May 2011, and I wish to express my thanks to the organisers of these events, especially Tom Leinster, Marino Gran and Enrico Vitale.  Readers familiar with this work may wish to note that it has not substantially changed since first being made available shortly after these conferences.


\ed
